\documentclass[journal]{IEEEtran}  

 \usepackage{etoolbox}
\makeatletter
\patchcmd{\@begintheorem}{\textit}{\textbf}{}{}
\makeatother
 \newtheorem{definition}{Definition}
 \newtheorem{problem}{Problem}
  \newtheorem{thm}{Theorem}
  
    \newtheorem{corollary}{Corollary}
  \newtheorem{remark}{Remark}
 \newtheorem{lemma}{Lemma}
  \newtheorem{prop}{Proposition}
\usepackage[]{amsmath}
\usepackage{amsfonts}
\usepackage{amssymb}
\usepackage{bbm}
\usepackage{breqn}
\usepackage{subcaption}
\usepackage{cite}
\usepackage{enumitem}
\usepackage{cases}
\usepackage{tabu}
\usepackage{url}
\usepackage{multicol}
\usepackage{color}
\usepackage[bottom]{footmisc}
\usepackage{environ}
\usepackage{tikz} 
\usetikzlibrary{automata,arrows,positioning,calc}
\usepackage[]{algpseudocode}
\algtext*{EndWhile}
\algtext*{EndIf}
\algtext*{EndFor}
\usepackage{algorithm}
\title{
Entropy Maximization for Markov Decision Processes Under Temporal Logic Constraints}
\author{Yagiz Savas, Melkior Ornik, Murat Cubuktepe, Mustafa O. Karabag, Ufuk Topcu  \thanks{This work was supported in part by the grants DARPA \# D19AP00004 and
AFRL \# FA9550-19-1-0169. } \thanks{ Y. Savas, M. Cubuktepe and U. Topcu are with the Department of Aerospace Engineering, University of Texas at Austin, TX, USA. E-mail: \{yagiz.savas, mcubuktepe, utopcu\}@utexas.edu   } \thanks{M. Ornik is with the Department of Aerospace Engineering and the Coordinated Science Laboratory, University of Illinois at Urbana-Champaign, IL, USA. E-mail: mornik@illinois.edu}  
 \thanks{M. O. Karabag is with the Department of Electrical and Computer Engineering, University of Texas at Austin, TX, USA. E-mail: karabag@utexas.edu  }}

\begin{document}

\maketitle
{\abstract  We study the problem of synthesizing a policy that maximizes the entropy of a Markov decision process (MDP) subject to a temporal logic constraint. Such a policy minimizes the predictability of the paths it generates, or dually, maximizes the exploration of different paths in an MDP while ensuring the satisfaction of a temporal logic specification. We first show that the maximum entropy of an MDP can be finite, infinite or unbounded. We provide necessary and sufficient conditions under which the maximum entropy of an MDP is finite, infinite or unbounded. We then present an algorithm which is based on a convex optimization problem to synthesize a policy that maximizes the entropy of an MDP. We also show that maximizing the entropy of an MDP is equivalent to maximizing the entropy of the paths that reach a certain set of states in the MDP. Finally, we extend the algorithm to an MDP subject to a temporal logic specification. In numerical examples, we demonstrate the proposed method on different motion planning scenarios and illustrate the relation between the restrictions imposed on the paths by a specification, the maximum entropy, and the predictability of paths.

  }
 
\section{Introduction}

Markov decision processes (MDPs) model sequential decision-making in stochastic systems with nondeterministic choices. A policy, i.e., a decision strategy, resolves the nondeterminism in an MDP and induces a stochastic process. In this regard, an MDP represents a (infinite) family of stochastic processes. In this paper, for a given MDP, we aim to synthesize a policy that induces a process with maximum entropy among the ones whose paths satisfy a temporal logic specification. 

Entropy, as an information-theoretic quantity, measures the unpredictability of outcomes in a random variable \cite{Cover}. Considering a stochastic process as an infinite sequence of (dependent) random variables, we define the entropy of a stochastic process as the joint entropy of these random variables by following \cite{Biondi},\cite{Chen}. Therefore, intuitively, our objective is to obtain a process whose paths satisfy a temporal logic specification in the most unpredictable way to an observer. 

Typically, in an MDP, a decision-maker is interested in satisfying certain properties \cite{Behcet} or accomplishing a task \cite{Jie}. Linear temporal logic (LTL) is a formal specification language \cite{Pneuli} that has been widely used to check the reliability of software \cite{Tan}, describe tasks for autonomous robots \cite{Gazit, Belta} and verify the correctness of communication protocols \cite{Kumar}. For example, in a robot navigation scenario, it allows to specify tasks such as safety (never visit the region A), liveness (eventually visit the region A) and priority (first visit the region A, then B). 

The entropy of paths of a (Markovian) stochastic process is introduced in \cite{Ekroot} and quantifies the randomness of realizations with fixed initial and final states. We first extend the definition for the entropy of paths to realizations that reach a certain set of states, rather than a fixed final state. Then, we show that the entropy of a stochastic process is equal to the entropy of paths of the process, if the process has a finite entropy. The established relation provides a mathematical basis to the intuitive idea that maximizing the entropy of an MDP minimizes the predictability of paths.

We observe that the maximum entropy of an MDP under stationary policies may not exist, i.e., for any given level of entropy, using stationary policies, one can induce a process whose entropy is greater than that level. In this case, we say that the maximum entropy of the MDP is unbounded. Additionally, if there exists a process with the maximum entropy, the entropy of such a process can be finite or infinite. Hence, before attempting to synthesize a policy that maximizes the entropy of an MDP, we first verify whether there exists a policy that attains the maximum entropy. 

The contributions of this paper are fourfold. First, we provide necessary and sufficient conditions on the structure of the MDP under which the maximum entropy of the MDP is finite, infinite or unbounded. We also present a polynomial-time algorithm to check whether the maximum entropy of an MDP is finite, infinite or unbounded. Second, we present a polynomial-time algorithm based on a convex optimization problem to synthesize a policy that maximizes the entropy of an MDP. Third, we show that maximizing the entropy of an MDP with non-infinite maximum entropy is equivalent to maximizing the entropy of paths of the MDP. Lastly, we provide a procedure to obtain a policy that maximizes the entropy of an MDP subject to a general LTL specification. 

The applications of this theoretical framework range from motion planning and stochastic traffic assignments to software security. In a motion planning scenario, for security purposes, an autonomous robot might need to randomize its paths while carrying out a mission \cite{Paruchuri, Paruchuri2}. In such a scenario, a policy synthesized by the proposed methods both provides probabilistic guarantees on the completion of the mission and minimizes the predictability of the robot's paths  through the use of online randomization mechanisms. Additionally, such a policy allows the robot to explore different parts of the environment \cite{Saerens}, and behave robustly against uncertainties in the environment \cite{Deep_learning}. The proposed methods can also be used to distribute traffic assignments over a network, which is known as stochastic traffic assignments \cite{Akamatsu}, as it promotes the use of different paths. Finally, as it is shown in \cite{Biondi}, the maximum information that an adversary can leak from a (deterministic) software, which is modeled as an MDP, can be quantified by computing the maximum entropy of the MDP. 

\textbf{Related Work.}  A preliminary version \cite{Yagiz} of this paper considered entropy maximization problem for MDPs subject to expected reward constraints. This considerably extended version includes an additional section establishing the relation between the maximum entropy of an MDP and the entropy of paths of the MDP, detailed proofs for all theoretical results, and additional numerical examples.

The computation of the maximum entropy of an MDP is first considered in \cite{Chen}, where the authors present a robust optimization problem to compute the maximum entropy for an MDP with finite maximum entropy. However, their approach does not allow to incorporate additional constraints due to the formulation of the problem. References \cite{Biondi} and \cite{Biondi2} compute the maximum entropy of an MDP for special cases without providing a general algorithm. 

The work \cite{Biondi} provides the necessary and sufficient conditions for an interval Markov chain (MC) to have a finite maximum entropy. Therefore, some of the results provided in this paper, e.g., the necessary and sufficient conditions for an MDP to have finite, unbounded or infinite maximum entropy, can be seen as an extension of the results given in \cite{Biondi}.

In \cite{Bullo,Bullo2}, the authors study the problem of synthesizing a transition matrix with maximum entropy for an irreducible MC subject to graph constraints. The problem studied in this paper is considerably different from that problem since MDPs represent a more general model than MCs, and an MC induced from an MDP by a policy is not necessarily irreducible.

In \cite{Paruchuri}, the authors maximize the entropy of a \textit{policy} while keeping the expected total reward above a threshold. They claim that the entropy maximization problem is not convex. Their formulation is a special case of the convex optimization problem that we provide in this paper. Therefore, here, we also prove the convexity of their formulation. 

The entropy of paths of absorbing MCs is discussed in \cite{Ekroot}, \cite{Akamatsu}, \cite{Kafsi}. The reference \cite{Saerens} establishes the equivalence between the entropy of paths and the entropy of an absorbing MC. We establish this relation for a general MC and show the connections to the maximum entropy of an MDP.

We also note that none of the above work discusses the unbounded and infinite maximum entropy for an MDP or considers LTL to specify desired system properties.

\textbf{Organization.} We provide the preliminary definitions and formal problem statement in Sections \ref{Prelim} and \ref{problem_set}, respectively. We analyze the properties of the maximum entropy of an MDP and present an algorithm to synthesize a policy that maximizes the entropy of an MDP in Section \ref{max-ent}. The relation between the maximum entropy of an MDP and the entropy of paths is established in Section \ref{relate_paths}. We present a procedure to synthesize a policy that maximizes the entropy of an MDP subject to an LTL specification in Section \ref{cons_section}. We provide numerical examples in Section \ref{examples_section} and conclude with suggestions for future work in Section \ref{conclusion_section}. Proofs for all results are provided in Appendix \ref{proofs_appendix}, and a procedure to synthesize a policy that maximizes the entropy of an MDP with infinite maximum entropy is presented in Appendix \ref{infinite_case_appendix}.

\section{Preliminaries}\label{Prelim}\noindent
\textbf{Notation:} For a set $S$,  we denote its power set and cardinality by $2^S$ and $\lvert S \rvert$, respectively. For a matrix $P$$\in$$\mathbb{R}^{n \times n}$, we use $P^k$ and $P^k_{i,j}$ to denote the k-th power of $P$ and the $(i,j)$-th component of the k-th power of $P$, respectively. All logarithms are to the base 2 and the  set $\mathbb{N}$ denotes $\{0,1,2,\ldots\}$.
\subsection{Markov chains and Markov decision processes}\label{MDP_subsection}
{\setlength{\parindent}{0cm}
\begin{definition}
A \textit{Markov decision process} (MDP) is a tuple $\mathcal{M}$$=$$(S, s_0, \mathcal{A}, \mathbb{P},\mathcal{AP},\mathcal{L})$ where $S$ is a finite set of states, $s_0$$\in$$S$ is the initial state, $\mathcal{A}$ is a finite set of actions, $ \mathbb{P}$$:$$S$$\times$$ \mathcal{A}$$ \times$$ S$$\rightarrow$$[0,1]$ is a transition function such that $\sum_{t\in S}\mathbb{P}(s,a,t)$$=$$1$ for all $s$$\in$$S$ and $a$$\in$$\mathcal{A}$, $\mathcal{AP}$ is a set of atomic propositions, and $\mathcal{L}$ $:$ $ S$$\rightarrow$$ 2^{\mathcal{AP}}$ is a function that labels each state with a subset of atomic propositions.
\end{definition}}
We denote the transition probability $ \mathbb{P}(s,a,t)$ by $ \mathbb{P}_{s,a,t}$, and all available actions in a state $s$$\in$$S$ by $\mathcal{A}(s)$. The set of successor states for a state action pair $(s,a)$ is defined as $Succ(s,a)$$:=$$\{t $$\in$$ S |  \mathbb{P}_{s,a,t}$$>$$0, a$$\in$$\mathcal{A}(s)\}$. The \textit{size of an MDP} is the number of triples $(s,a,t)$$\in$$S$$\times$$\mathcal{A}$$\times$$S$ such that $\mathbb{P}_{s,a,t}$$>$$0$. 

A \textit{Markov chain} (MC) $\mathcal{C}$ is an MDP such that $\lvert \mathcal{A}\rvert$$=$$1$. We denote the transition function (matrix) for an MC by $\mathcal{P}$, and the set of successor states for a state $s$$\in$$S$ by $Succ(s)$$=$$\{t $$\in$$ S | \mathcal{P}_{s,t}$$>$$0\}$. The \textit{expected residence time} in a state $s$$\in$$ S$ for an MC $\mathcal{C}$ is defined as
\begin{align}
\label{residence}
\xi_s:=\sum_{k=0}^{\infty}\mathcal{P}^k_{s_0,s}.\end{align}
The expected residence time $\xi_s$ represents the expected number of visits to state $s$ starting from the initial state ~\cite{Marta}. A state $s$$\in$$S$ is \textit{recurrent} for an MC if and only if $\xi_s$$=$$\infty$, and is \textit{transient} otherwise; it is \textit{stochastic} if and only if it satisfies $\lvert Succ(s)\rvert$$>$$1$, and is \textit{deterministic} otherwise; and it is \textit{reachable} if and only if $\xi_s$$>$$0$, and is \textit{unreachable} otherwise.  
{\setlength{\parindent}{0cm}
\noindent \begin{definition}
A \textit{policy} for an MDP $\mathcal{M}$ is a sequence $\pi$$=$$\{\mu_0, \mu_1, \ldots\}$ where each $\mu_k $$:$$ S $$ \times$$ \mathcal{A}$$\rightarrow$$[0,1]$ is a function such that $\sum_{a\in \mathcal{A}(s)}\mu_k(s,a)$$=$$1$ for all $s$$\in$$S$. A \textit{stationary} policy is a policy of the form $\pi$$=$$\{\mu, \mu, \ldots\}$. For an MDP $\mathcal{M}$, we denote the set of all policies and all stationary policies by $\Pi(\mathcal{M})$ and $\Pi^S(\mathcal{M})$, respectively.
\end{definition}}
 We denote the probability of choosing an action $a$$\in$$\mathcal{A}(s)$ in a state $s$$\in$$S$ under a stationary policy $\pi$ by $\pi_s(a)$. For an MDP $\mathcal{M}$, a stationary policy $\pi$$\in$$\Pi^S(\mathcal{M})$ induces an MC denoted by $\mathcal{M}^{\pi}$. We refer to $\mathcal{M}^{\pi}$ as \textit{induced MC} and specify the transition matrix for $\mathcal{M}^{\pi}$ by $\mathcal{P}^{\pi}$, whose $(s,t)$-th component is given by
\begin{align}
\label{induced_MC_prob}
\mathcal{P}^{\pi}_{s,t}=\sum_{a\in\mathcal{A}(s)}\pi_s(a) \mathbb{P}_{s,a,t}.
\end{align}
Throughout the paper, we assume that for a given MDP $\mathcal{M}$, for any state $s$$\in$$S$ there exists an induced MC $\mathcal{M}^{\pi}$ for which the state $s$ is reachable. This is a standard assumption for MDPs \cite{Belta}, which ensures that each state in the MDP is reachable under some policy.

An infinite sequence $\varrho^{\pi}$$=$$s_0s_1s_2\ldots$ of states generated in $\mathcal{M}$ under a policy $\pi$$\in$$\Pi(\mathcal{M})$ is called a \textit{path}, starting from the initial state $s_0$ and satisfies $\sum_{a_k\in \mathcal{A}(s_k)}\mu_k(s_k)(a_k)\mathbb{P}_{s_k,a_k,s_{k+1}}$$>$$0$ for all $k$$\geq$$0$. Any finite prefix of $\varrho^{\pi}$ that ends in a state is a finite path fragment. We define the set of all paths and finite path fragments in $\mathcal{M}$ under the policy $\pi$ by $Paths^{\pi}(\mathcal{M})$ and $Paths_{fin}^{\pi}(\mathcal{M})$, respectively. 

We use the standard probability measure over the outcome set $Paths^{\pi}(\mathcal{M})$ \cite{Model_checking}. For a path $\varrho^{\pi}$$\in$$Paths^{\pi}(\mathcal{M})$, let the sequence $s_0s_1\ldots s_n$ be the finite path fragment of length $n$, and let $Paths^{\pi}(\mathcal{M})(s_0s_1\ldots s_n)$ denote the set of all paths in $Paths^{\pi}(\mathcal{M})$ starting with the prefix $s_0s_1\ldots s_n$. The probability measure $\text{Pr}_{\mathcal{M}}^{\pi}$ defined on the smallest $\sigma$-algebra over $Paths^{\pi}(\mathcal{M})$ that contains $Paths^{\pi}(\mathcal{M})(s_0s_1\ldots s_n)$ for all $s_0s_1\ldots s_n$$\in$$Paths_{fin}^{\pi}(\mathcal{M})$ is the unique measure that satisfies 
{\setlength{\mathindent}{0pt}
\begin{flalign}
\label{first_measure}\noindent
\text{Pr}_{\mathcal{M}}^{\pi}\{Paths^{\pi}&(\mathcal{M})(s_0\ldots s_n)\}=\nonumber \\
&\prod_{0\leq k < n} \sum_{a_k\in\mathcal{A}(s_k)}\mu_k(s_k)(a_k) \mathbb{P}_{s_k,a_k, s_{k+1}}.
\end{flalign}}\noindent\vspace{-0.7cm}

\subsection{The entropy of stochastic processes}\label{entropy_subsection}
For a (discrete) random variable $X$, its support $\mathcal{X}$ defines a countable sample space from which $X$ takes a value $x$$\in$$\mathcal{X}$ according to a probability mass function (pmf) $p(x)$$:=$$\text{Pr}(X$$=$$x)$. The \textit{entropy} of a random variable $X$ with countable support $\mathcal{X}$ and pmf $p(x)$ is defined as
\begin{align}\centering
H(X):=-\sum_{x\in\mathcal{X}} p(x)\log p(x).
\end{align}
We use the convention that $0$$\log$$0$$=$$0$. Let $(X_0,X_1)$ be a pair of random variables with the joint pmf $p(x_0,x_{1})$ and the support $\mathcal{X}\times \mathcal{X}$. The \textit{joint entropy} of $(X_0,X_1)$ is 
\begin{align}
\label{joint_entropy}
H(X_0,X_1):= -\sum_{x_0\in \mathcal{X}}\sum_{x_{1}\in \mathcal{X}}p(x_0,x_{1})\log p(x_0,x_{1}),
\end{align}\noindent
and the \textit{conditional entropy} of $X_1$ given $X_0$ is
\begin{align}
\label{conditional_entropy}
&H(X_1 | X_0):=-\sum_{x_0\in \mathcal{X}}\sum_{x_{1}\in \mathcal{X}}p(x_0,x_{1})\log p(x_1 |x_0).
\end{align}\noindent
The definitions of the joint and conditional entropies extend to collection of $k$ random variables as it is shown in \cite{Cover}. 
A discrete \textit{stochastic process} $\mathbb{X}$ is a discrete time-indexed sequence of random variables, i.e., $\mathbb{X}$$=$$\{X_k$$\in$$\mathcal{X}$ $:$ $k$$\in$$\mathbb{N}\}$. 
{\setlength{\parindent}{0cm}
\noindent\begin{definition} (Entropy of a stochastic process) \cite{Biondi_thesis}
The \textit{entropy of a stochastic process} $\mathbb{X}$ is defined as 
 \begin{align}\label{entropy_def_stochastic}
 H(\mathbb{X}):=\lim_{k\rightarrow \infty}H( X_0, X_1\ldots, X_{k}).
 \end{align}
 \end{definition}}
Note that this definition is different from the \textit{entropy rate} of a stochastic process, which is defined as $\lim_{k\rightarrow \infty}\frac{1}{k}H( X_0, X_1\ldots, X_{k})$ when the limit exists \cite{Cover}. The limit in \eqref{entropy_def_stochastic} either converges to a non-negative real number or diverges to positive infinity \cite{Biondi_thesis}. 

An MC $\mathcal{C}$ is equipped with a discrete stochastic process $\{X_k$$\in$$S$ $:$ $ k$$\in$$\mathbb{N}\}$ where each $X_k$ is a random variable over the state space $S$. For a given k-dimensional pmf $p(s_0,s_1,\ldots, s_k)$, this process respects the \textit{Markov property}, i.e., $p(s_k|s_{k-1},\ldots,s_0)$$=$$p(s_k|s_{k-1})$
for all $k$$\in$$\mathbb{N}$. Then, the \textit{entropy of a Markov chain} $\mathcal{C}$ is given by 
\begin{align}
\label{process_entropy}
H(\mathcal{C})=H(X_0)+\sum_{i=1}^{\infty}H(X_{i}| X_{i-1})
\end{align} 
using \eqref{joint_entropy}, \eqref{conditional_entropy} and \eqref{entropy_def_stochastic}. Note that $H(X_0)$$=$$0$, since we define an MC with a unique initial state.

For an MDP $\mathcal{M}$, a policy $\pi$$\in$$\Pi(\mathcal{M})$ induces a discrete stochastic process $\{X_k $$\in$$S$ $:$ $ k$$\in$$\mathbb{N}\}$.  We denote the entropy of an MDP $\mathcal{M}$ under a policy $\pi$$\in$$\Pi(\mathcal{M})$ by $H(\mathcal{M},\pi)$. 
Using the next proposition, we restrict our attention to stationary policies for maximizing the entropy of an MDP. 
{\setlength{\parindent}{0cm}\noindent
\begin{prop}
\label{memoryless_1}
The following equality holds:
\vspace{-0.2cm}
\begin{align}\label{eqeq}
\sup_{\pi\in\Pi(\mathcal{M})}H(\mathcal{M},{\pi})=\sup_{\pi\in\Pi^S(\mathcal{M})}H(\mathcal{M},{\pi}).\quad \triangleleft
\end{align}
\end{prop}}A proof for Proposition \ref{memoryless_1} is provided in Appendix \ref{proofs_appendix}.
{\setlength{\parindent}{0cm}\noindent\begin{remark}  If the supremum in \eqref{eqeq} is infinite, the set of stationary policies may not be sufficient to attain the supremum while a non-stationary policy can attain it. In particular, there exists a family of distributions that are defined over a countable support and have infinite entropy (see equation (7) in \cite{Bacetti}). It can be shown that for some MDPs, there exists a non-stationary policy that induces a stochastic process with such a probability distribution, and hence, have infinite entropy, while stationary policies can only induce stochastic processes with finite entropies\footnote{ A preliminary version \cite{Yagiz} of this paper relied on Proposition 36 from \cite{Chen}. This proposition is not valid in general. Here, we provide the corrected results by defining the maximum entropy of an MDP over stationary policies.}.\end{remark}}\vspace{-0.2cm}
{\setlength{\parindent}{0cm}\noindent\begin{definition} (Maximum entropy of an MDP) The \textit{maximum entropy of an MDP} $\mathcal{M}$ is
\begin{align}\label{max_ent_definition}
H(\mathcal{M}):=\sup_{\pi\in\Pi^{S}(\mathcal{M})}H(\mathcal{M},{\pi}).
\end{align}\end{definition}}
A policy $\pi^{\star}$$\in$$\Pi^S(\mathcal{M})$ \textit{maximizes} the entropy of an MDP $\mathcal{M}$ if $H(\mathcal{M})$$=$$H(\mathcal{M},{\pi^{\star}})$. Finally, we define the properties of the maximum entropy of an MDP as follows.{\setlength{\parindent}{0cm}
\noindent \begin{definition}(The properties of the maximum entropy)\label{props_def}
The maximum entropy of an MDP $\mathcal{M}$ is
\begin{itemize} 
\item \textit{finite}, if and only if
\begin{align}\label{def_finite_ent}
H(\mathcal{M})=\max_{\pi\in\Pi^{S}(\mathcal{M})}H(\mathcal{M},{\pi})<\infty;
\end{align}
\item \textit{infinite}, if and only if
\begin{align} H(\mathcal{M})=\max_{\pi\in\Pi^{S}(\mathcal{M})}H(\mathcal{M},{\pi})=\infty; 
\end{align}
\item \textit{unbounded}, if and only if the following two conditions hold.
\begin{align}\label{unbounded_definition_1}
(i)\qquad &H(\mathcal{M})=\sup_{\pi\in\Pi^{S}(\mathcal{M})}H(\mathcal{M},{\pi})=\infty,\\ \label{unbounded_definition_2}
(ii)\qquad &H(\mathcal{M},{\pi})<\infty\ \  \text{for all}\ \pi\in\Pi^{S}(\mathcal{M}).
\end{align}
\end{itemize}
\end{definition}}
Although it is not defined here, there is a fourth possible property which is unachievable finite maximum entropy, i.e., $\max_{\pi\in\Pi^S(\mathcal{M})}H(\mathcal{M},{\pi})$$<$$H(\mathcal{M})$$<$$\infty$. In Theorem \ref{Theorem1}, we show that it is not possible for the maximum entropy of an MDP to have this property.
\subsection{Linear temporal logic}
We employ linear temporal logic (LTL) to specify tasks and refer the reader to \cite{Model_checking} for the syntax and semantics of LTL. 

An LTL formula is built up from a set $\mathcal{AP}$ of atomic propositions, logical connectives such as conjunction ($\land$) and negation ($\lnot$), and temporal modal operators such as always ($\square$) and eventually ($\lozenge$). An infinite sequence of subsets of $\mathcal{AP}$ defines an infinite \textit{word}, and an LTL formula is interpreted over infinite words on $2^{\mathcal{AP}}$. We denote by $w$$\models$$\varphi$ that a word $w$$=$$w_0w_1w_2\ldots$ satisfies an LTL formula $\varphi$.
{\setlength{\parindent}{0cm}
\noindent \begin{definition}
A \textit{deterministic Rabin automaton} (DRA) is a tuple $A$$=$$(Q,q_0, \Sigma, \delta, Acc)$ where $Q$ is a finite set of states, $q_0$$\in$$Q$ is the initial state, $\Sigma$ is the alphabet, $\delta$$:$$Q$$\times$$\Sigma$$\rightarrow$$Q$ is the transition relation, and $Acc$$\subseteq$$2^{Q}$$\times$$2^{Q}$ is the set of accepting state pairs. 
\end{definition}}
A \textit{run} of a DRA $A$, denoted by $\sigma$$=$$q_0q_1q_2\ldots$, is an infinite sequence of states in $\mathcal{A}$ such that for each $i$$\geq$$0$, $q_{i+1}$$\in$$\delta(q_i, p)$ for some $p$$\in$$\Sigma$. A run $\sigma$ is \textit{accepting} if there exists a pair $(J,K)$$\in$$Acc$ and an $n$$\geq$$0$ such that (i) for all $m$$\geq$$n$ we have $q_m$$\not\in$$J$, and (ii) there exists infinitely many $k$ such that $q_k$$\in$$K$.

For any LTL formula $\varphi$ built up from $\mathcal{AP}$, a DRA $A_{\varphi}$ can be constructed with input alphabet $2^{\mathcal{AP}}$ that accepts all and only words over $\mathcal{AP}$ that satisfy $\varphi$ \cite{Model_checking}. 

For an MDP $\mathcal{M}$ under a policy $\pi$, a path $\varrho^{\pi}$$=$$s_0s_1\ldots$ generates a word $w$$=$$w_0w_1\ldots$ where $w_k$$=$$\mathcal{L}(s_k)$ for all $k$$\geq$$0$. With a slight abuse of notation, we use $\mathcal{L}(\varrho^{\pi})$ to denote the word generated by $\varrho^{\pi}$. For an LTL formula $\varphi$, the set $\{\varrho^{\pi}$$\in$$ Paths^{\pi}(\mathcal{M})$$:$$ \mathcal{L}(\varrho^{\pi})$$\models$$\varphi\}$ is measurable \cite{Model_checking}. We define 
\begin{equation*}
\begin{aligned}
\text{Pr}_{\mathcal{M}}^{\pi}(s_0\models \varphi):=\text{Pr}_{\mathcal{M}}^{\pi}\{\varrho^{\pi}\in Paths^{\pi}(\mathcal{M}) : \mathcal{L}(\varrho^{\pi})\models \varphi\}
\end{aligned}
\end{equation*}
as the probability of satisfying the LTL formula $\varphi$ for an MDP $\mathcal{M}$ under the policy $\pi$$\in$$\Pi(\mathcal{M})$.  
\section{Problem Statement}\label{problem_set}
The first problem we study concerns the synthesis of a policy that maximizes the entropy of an MDP.
 {\setlength{\parindent}{0cm}
\noindent \begin{problem}(\textbf{Entropy Maximization})\label{prob_1}
For a given MDP $\mathcal{M}$, provide an algorithm to verify whether there exists a policy $\pi^{\star}$$\in$$\Pi^S(\mathcal{M})$ such that $H(\mathcal{M})$$=$$H(\mathcal{M},\pi^{\star})$. If such a policy exists, provide an algorithm to synthesize it. If it does not exist, provide a procedure to synthesize a policy $\pi'$$\in$$\Pi^S(\mathcal{M})$ such that $H(\mathcal{M},\pi')$$\geq$$\ell$ for a given constant $\ell$. 
\end{problem}}

For an MDP $\mathcal{M}$, the synthesis of a policy $\pi'$$\in$$\Pi^S(\mathcal{M})$ such that $H(\mathcal{M},\pi')$$\geq$$\ell$ allows one to induce a stochastic process with the desired level of entropy, even if there exists no stationary policy that maximizes the entropy of $\mathcal{M}$.

In the second problem, we introduce linear temporal logic (LTL) specifications to the framework. In particular, we consider the problem of synthesizing a policy that induces a stochastic process with maximum entropy whose paths satisfy a given LTL formula with desired probability. The formal statement of the second problem is deferred to Section \ref{cons_section} since it requires the introduction of additional notations.

\section{Entropy maximization for MDPs }
\label{max-ent}
 In this section, we focus on the entropy maximization problem. We refer to a policy as an \textit{optimal} policy for an MDP if it maximizes the entropy of the MDP.
 \subsection{The entropy of MCs versus MDPs}
 For an MC, the \textit{local entropy} of a state $s$$\in$$S$ is defined as
\begin{align}
\label{local_entropy_def}
L(s):=-\sum_{t\in S}\mathcal{P}_{s,t}\log \mathcal{P}_{s,t}.
\end{align}
The following proposition characterizes the relationship between the local entropy of states and the entropy of an MC.{\setlength{\parindent}{0cm}
\noindent \begin{prop}
\label{Biondi_theorem}
(Theorem 1 in \cite{Biondi}) For an MC $\mathcal{C}$,
\begin{align} 
\label{MC_entropy_biondi}
H(\mathcal{C})=\sum_{s\in S}L(s)\xi_s.\quad \triangleleft
\end{align}
\end{prop}}\noindent
\par An MC $\mathcal{C}$ has a finite entropy if and only if all of its recurrent states have zero local entropy \cite{Biondi}. That is, $H(\mathcal{C})$$<$$\infty$ if and only if for all states $s$$\in$$ S$, $\xi_s$$=$$\infty$ implies $L(s)$$=$$0$. If the entropy of an MC is finite, each recurrent state $s$$\in$$S$ has a single successor state, i.e., $\lvert Succ(s)\rvert$$=$$1$. Consequently, recurrent states have no contribution to the sum in \eqref{process_entropy}. In this case, we take the sum in \eqref{MC_entropy_biondi} only over the transient states.

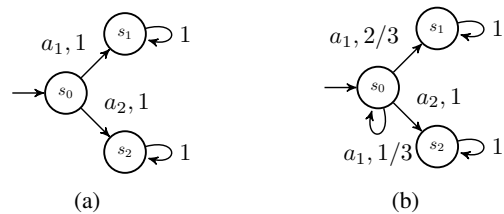
\begin{figure}[b]\vspace{-0.4cm}
\hspace{0.05\linewidth}
\begin{subfigure}[b]{0.3\linewidth}
\centering
\scalebox{0.9}{
\begin{tikzpicture}[->, >=stealth', auto, semithick, node distance=1.5cm]

    \tikzstyle{every state}=[fill=white,draw=black,thick,text=black,scale=0.7]

    \node[state,initial,initial text=] (s_0) {$s_0$};
    \node[state] (s_1) [above right =6mm of s_0]  {$s_1$};
    \node[state] (s_2) [below right=6mm of s_0]  {$s_2$};

\path
(s_0)  edge     node{$a_1, 1$}     (s_1)
(s_0)	 edge      node{$a_2, 1$}     (s_2)
(s_1)	 edge    [loop right]       node{$1$}     (s_1)
(s_2)	 edge    [loop right]       node{$ 1$}     (s_2);
\end{tikzpicture}}
\caption{}
\label{fig:MDP_1}
\end{subfigure}
\hspace{0.1\linewidth}
\begin{subfigure}[b]{0.4\linewidth}
\centering
\scalebox{0.9}{
\begin{tikzpicture}[->, >=stealth', auto, semithick, node distance=2cm]

    \tikzstyle{every state}=[fill=white,draw=black,thick,text=black,scale=0.7]

    \node[state,initial,initial text=] (s_0) {$s_0$};
    \node[state] (s_1) [above right =6mm of s_0]  {$s_1$};
    \node[state] (s_2) [below right=6mm of s_0]  {$s_2$};

\path
(s_0)  edge      [loop below] node{$a_1, 1/3$}     (s_0)
(s_0)  edge     node{$a_1, 2/3$}     (s_1)
(s_0)	 edge      node{$a_2, 1$}     (s_2)
(s_1)	 edge    [loop right]       node{$1$}     (s_1)
(s_2)	 edge    [loop right]       node{$ 1$}     (s_2);
\end{tikzpicture}}
\caption{}
\label{fig:MDP_2}
\end{subfigure}
\caption{Randomizing actions uniformly at each state may or may not achieve the maximum entropy. The optimal policy for the MDP given in (a) is $\pi_{s_0}(a_1)$$=$$\pi_{s_0}(a_2)$$=$$1/2$, and for the MDP given in (b) is $\pi_{s_0}(a_1)$$=$$2/3$, $\pi_{s_0}(a_2)$$=$$1/3$.}
\label{fig:examples_policies}
\end{figure}

For an MDP, different policies may induce stochastic processes with different entropies. For example, consider the MDP given in Fig. \ref{fig:MDP_1} and suppose that the action $a_1$ at state $s_0$ is taken with probability $\varepsilon$. If we let $\varepsilon$ range over $[0, \frac{1}{2}]$, then the entropy of the resulting stochastic processes ranges over $[0,1]$. The optimal policy for this MDP is $\pi_{s_0}(a_1)$$=$$\pi_{s_0}(a_2)$$=$$1/2$, which uniformly randomizes actions. 
\par Unlike the MDP given in Fig. \ref{fig:MDP_1}, the maximum entropy of an MDP is not generally achieved by a policy that chooses available actions at each state uniformly. For example, consider the MDP given in Fig. \ref{fig:MDP_2}. The optimal policy for this MDP is $\pi_{s_0}(a_1)$$=$$2/3$, $\pi_{s_0}(a_2)$$=$$1/3$.

\par Examples given in Fig. \ref{fig:examples_policies} show that finding an optimal policy for an MDP may not be trivial. To analyze the maximum entropy of an MDP, we first obtain a compact representation of the maximum entropy as follows. For an MC $\mathcal{M}^{\pi}$ induced from an MDP $\mathcal{M}$ by a policy $\pi$$\in$$\Pi^S(\mathcal{M})$, let the expected residence time in a state $s$$\in$$S$ be 
\begin{align}
\label{MDP_residence_time}
\xi^{\pi}_s:=\sum_{k=0}^{\infty}(\mathcal{P}^{\pi})^k_{s_0,s}.
\end{align}

Additionally, let the local entropy of a state $s$$\in$$S$ in $\mathcal{M}^{\pi}$ be
$L^{\pi}(s)$$:=$$-\sum_{t\in S}\mathcal{P}^{\pi}_{s,t}\log\mathcal{P}^{\pi}_{s,t}$. Then, the maximum entropy of $\mathcal{M}$ can be written as
\begin{flalign}
\label{objective_MDP}
&H(\mathcal{M})=\sup_{\pi\in{\Pi^S(\mathcal{M})}}\Big[\sum_{s\in S}\xi^{\pi}_sL^{\pi}(s)\Big].
\end{flalign}
Note that the right hand side of \eqref{objective_MDP} can still be infinite or unbounded. We analyze the properties of the maximum entropy of MDPs in the next section.
\subsection{Properties of the maximum entropy of MDPs}\label{characteristics}

\begin{figure}[b!]\vspace{-0.5cm}
\begin{subfigure}[t]{0.33\linewidth}
\scalebox{0.9}{
\begin{tikzpicture}[->, >=stealth', auto, semithick, node distance=2cm]

    \tikzstyle{every state}=[fill=white,draw=black,thick,text=black,scale=0.7]

    \node[state,initial,initial text=] (s_0) {$s_0$};
    \node[state] (s_1) [right=10mm of s_0]  {$s_1$};

\path
(s_0)  edge  [loop above=10]    node{$a_1, 1$}     (s_0)
(s_0)	 edge     node{$a_2, 1$}     (s_1)
(s_1)  edge  [loop right=10]    node{$1$}     (s_1);
\end{tikzpicture}}
\caption{}
\label{fig:unbounded}
\end{subfigure}
\hspace{0.13\linewidth}
\begin{subfigure}[t]{0.33\linewidth}
\scalebox{0.9}{
\begin{tikzpicture}[->, >=stealth', auto, semithick, node distance=2cm]

    \tikzstyle{every state}=[fill=white,draw=black,thick,text=black,scale=0.7]

    \node[state,initial,initial text=] (s_0) {$s_0$};
    \node[state] (s_1) [right=10mm of s_0]  {$s_1$};

\path
(s_0)  edge  [loop above=10]    node{$a_1, 1$}     (s_0)
(s_0)  edge  [bend left=15]    node{$a_2, 1$}     (s_1)
(s_1)  edge  [loop right=10]    node{$a_1,1$}     (s_1)
(s_1)  edge  [bend left=15]    node{$a_2,1$}     (s_0);
\end{tikzpicture}}
\caption{}
\label{fig:infinite}
\end{subfigure}
\caption{Examples of MDPs with (a) unbounded maximum entropy and (b) infinite maximum entropy.}
\label{fig:examples}
\end{figure}
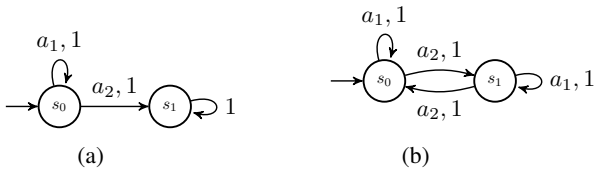

The maximum entropy of an MDP can be infinite or unbounded even for simple cases. For example, consider MDPs given in Fig. \ref{fig:examples}. For the MDP shown in Fig. \ref{fig:unbounded}, let the action $a_2$ be taken with probability $\delta$$\in$$(0,1]$ in state $s_0$. Then, the expected residence time $\xi^{\pi}_{s_0}$ in state $s_0$ is equal to $\frac{1}{\delta}$, and the entropy of the induced MC $\mathcal{M}^{\pi}$ is given by
\begin{align}
H(\mathcal{M},{\pi})=-\frac{(1-\delta)\log(1-\delta)+\delta\log(\delta)}{\delta},
\end{align} 
which satisfies $H(\mathcal{M},{\pi})$$\rightarrow$$\infty$ as $\delta$$\rightarrow$$0$. Note also that if $\delta$$=$$0$, the entropy of the induced MC is zero due to \eqref{MC_entropy_biondi}. Hence, the maximum entropy is unbounded, and there is no optimal stationary policy for this MDP. 

For the MDP given in Fig. \ref{fig:infinite}, choosing a policy such that $\pi_i(a_j)$$>$$0$ for $i$$=$$1,2$, $j$$=$$1,2$ yields $\xi^{\pi}_{s_0}$$=$$\xi^{\pi}_{s_1}$$=$$\infty$ and $L^{\pi}(s_0)$$>$$0$, $L^{\pi}(s_1)$$>$$0$. Then, the maximum entropy of this MDP is infinite, and the maximum can be attained by any randomized policy. 

\par Examples in Fig. \ref{fig:examples} show that we should first verify the existence of optimal policies before attempting to synthesize them. We need the following definitions about the structure of MDPs to state the conditions that cause an MDP to have finite, infinite or unbounded maximum entropy.

A \textit{directed graph} (digraph) is a tuple $G$$=$$(V,E)$ where $V$ is a set of vertices and $E$ is a set of ordered pairs of vertices $V$. For a digraph $G$, a path $v_1v_2\ldots v_n$ from vertex $v_1$ to $v_n$ is a sequence of vertices such that $(v_k,v_{k+1})$$\in$$E$ for all $1$$\leq$$k$$<$$n$. A digraph $G$ is \textit{strongly connected} if for every pair of vertices $u,v$$\in$$V$, there is a path from $u$ to $v$, and $v$ to $u$. 

A \textit{sub-MDP} of an MDP is a pair $(C,D)$ where $\emptyset$$\neq$$C$$\subseteq$$S$ and $D$$: $$C$$\rightarrow$$ 2^{\mathcal{A}}$ is a function such that (i) $D(s)$$\subseteq$$\mathcal{A}(s)$ is non-empty for all $s$$\in$$C$, and (ii) $s$$\in$$C$ and $a$$\in$$D(s)$ imply that $Succ(s,a)$$\subseteq$$C$. An \textit{end component} is a sub-MDP $(C,D)$ such that the digraph induced by $(C,D)$ is strongly connected.
{\setlength{\parindent}{0cm}
\noindent \begin{definition}
A \textit{maximal end component} (MEC) $(C,D)$ in an MDP is an end component such that there is no end component $(C', D')$ with $(C,D)$$\neq$$(C',D')$, and $C$$\subseteq$$C'$ and $D(s)$$\subseteq$$D'(s)$ for all $s$$\in$$C$.\end{definition}}
 \par A MEC $(C,D)$ in an MDP is \textit{bottom strongly connected} (BSC) if for all $s$$\in$$ C$, $\mathcal{A}(s)$$\backslash$$D(s)$$=$$\emptyset$. For a given state $s$$\in$$C$, we define the set of all actions under which the MDP can leave the MEC $(C, D)$ as $D_0(s)$$:=$$\{a$$\in$$\mathcal{A}(s) | Succ(s,a)$$\not\subseteq$$ C\}$. Note that in a BSC MEC $(C,D)$, $D_0(s)$$=$$\emptyset$ for all $s$$\in$$C$. 
  {\setlength{\parindent}{0cm}
 \noindent \begin{lemma}\label{Successor}
For an MDP $\mathcal{M}$ with MECs $(C_i,D_i)$ $i$$=$$1,2,\ldots,n$, let $C$$:=$$\cup_{i=1}^nC_i$ and $D$$:=$$\cup_{i=1}^nD_i$. Then, there exists an induced MC $\mathcal{M}^{\pi}$ for which a state $s$$\in$$C$ is both stochastic and recurrent if and only if $|\cup_{a\in D(s)}Succ(s,a)|$$>$$1$.$\quad\triangleleft$ \end{lemma}}
{\setlength{\parindent}{0cm}
 \begin{thm}\label{Theorem1}
For an MDP $\mathcal{M}$ with MECs $(C_i,D_i)$ $i$$=$$1,2,\ldots,n$, let $C$$:=$$\cup_{i=1}^nC_i$ and $D$$:=$$\cup_{i=1}^nD_i$. Then, the following statements hold. \\
(i)  $H(\mathcal{M})$  is infinite if and only if there exists an induced MC for which a state $s$$\in$$C$ is both stochastic and recurrent.\\
(ii)  $H(\mathcal{M})$ is unbounded if and only if $|\cup_{a\in D(s)}Succ(s,a)|$$=$$1$ for all $s$$\in$$C$, and there exists a MEC that is not bottom strongly connected. \\
(iii)  $H(\mathcal{M})$ is finite if and only if it is not infinite and not unbounded.$\quad \triangleleft$
\end{thm}}\noindent

Proofs for above results can be found in Appendix \ref{proofs_appendix}. Informally, Theorem \ref{Theorem1} states that for an MDP to have finite maximum entropy, all recurrent states of all MCs that are induced from the MDP by a stationary policy should be deterministic. Although necessary conditions for the finiteness of the maximum entropy is quite restrictive, there are some special cases, such as stochastic shortest path (SSP) problems  \cite{Bertsekas}, where MDP structures actually satisfy the necessary conditions. Specifically, since all \textit{proper} policies in SSP problems are guaranteed to reach an absorbing target state within finite time steps with probability 1, the problem of synthesizing a \text{proper} policy with maximum entropy has a finite solution.

The following corollary is due to Proposition \ref{memoryless_1}, Theorem \ref{Theorem1}, and the definition of finite maximum entropy \eqref{def_finite_ent}.
{\setlength{\parindent}{0cm}
 \begin{corollary}\label{corollary1} If $ \sup_{\pi\in\Pi(\mathcal{M})}H(\mathcal{M},{\pi})$$<$$\infty$, then we have
 \begin{align}
 \sup_{\pi\in\Pi(\mathcal{M})}H(\mathcal{M},{\pi})= \max_{\pi\in\Pi^S(\mathcal{M})}H(\mathcal{M},{\pi}).
 \end{align}
 \end{corollary}}

\begin{algorithm}[b]
 \caption{Verify the properties of the maximum entropy.}
 \textbf{Require:} $\mathcal{M}$$=$$(S,s_0, \mathcal{A},\mathbb{P},\mathcal{AP},\mathcal{L})$\\
  \textbf{Return:} R

Find: MECs $(C_i,D_i)$, $i=1,...,n$
  
Find: $Succ(s,a)$ for all $s$$\in$$S$, $a$$\in$$\mathcal{A}(s)$ 

R := $\emptyset$;

  \begin{algorithmic}
\For {$i$$=$$ 1,2,\ldots,n$}
 	\For {$s$ \textbf{in} $C_i$}
		\If{$\lvert \cup_{a\in D_i(s)} Succ(s,a)\rvert$$>$$1$}
			\State  R := R $\cup$ $\{\text{infinite}\}$ ;
		\EndIf
		 \If{ $\mathcal{A}(s)\backslash D_i(s)$$\neq$$\emptyset$ }
			\State  R := R $\cup$ $\{\text{unbounded}\}$ ;
		\EndIf

	\EndFor
 \EndFor 
 \If{\qquad \ \ \text{infinite}$\in$R } R=infinite 
 \ElsIf {\ \ \text{unbounded}$\in$R }\ R=unbounded
 \Else \qquad R=finite
 \EndIf
 \end{algorithmic}
 \label{algo_1}
\end{algorithm}

We present Algorithm \ref{algo_1} which, for an MDP $\mathcal{M}$, verifies whether $H(\mathcal{M})$ is finite, infinite or unbounded by checking the necessary conditions in Theorem \ref{Theorem1}. For $\mathcal{M}$, its MECs can be found in $\mathcal{O}(\lvert S\rvert ^2)$ time \cite{Model_checking}, $Succ(s,a)$ can be found in $\mathcal{O}(\lvert S\rvert ^2\lvert \mathcal{A}\rvert)$ time, and the necessary conditions can be verified in $\mathcal{O}(\lvert S\rvert)$ time since no state can belong to more than one MEC. Hence, Algorithm \ref{algo_1} runs in polynomial-time in the size of $\mathcal{M}$.

\subsection{Policy synthesis  }\label{policy_syntesis_section}
We now provide algorithms to synthesize policies that solve the entropy maximization problem. 

\subsubsection{Finite maximum entropy}\label{finite_policy_synthesis}
We first modify a given MDP by making all states in its MECs absorbing.
{\setlength{\parindent}{0cm}
\begin{prop}\label{modified_MDP_prop} Let $\mathcal{M}$ be an MDP such that $H(\mathcal{M})$$<$$\infty$, $(C_i,D_i)$ $i$$=$$1,2,\ldots,n$ be MECs in $\mathcal{M}$, $C$$:=$$\cup_{i=1}^nC_i$, and $\mathcal{M'}$ be the modified MDP that is obtained from $\mathcal{M}$ by making all states $s$$\in$$C$ absorbing, i.e., if $s$$\in$$C$, then $\mathbb{P}_{s,a,s}$$=$$1$ for all $a$$\in$$\mathcal{A}(s)$ in $\mathcal{M'}$. Then, we have
$ H(\mathcal{M})$$=$$H(\mathcal{M'})$.$\quad \triangleleft$ 
\end{prop}}

There is a one-to-one correspondence between the paths of $\mathcal{M}$ and $\mathcal{M'}$ since all states in the set $C$ must have a single successor state in an MDP with finite maximum entropy due to Theorem \ref{Theorem1}. Moreover, for a given policy $\pi'$$\in$$\Pi^S(\mathcal{M'})$ on $\mathcal{M'}$, the policy $\pi$$\in$$\Pi^S(\mathcal{M})$ induced by $\pi'$ on $\mathcal{M}$ is the same policy with $\pi'$, i.e. $\pi'$$=$$\pi$. Therefore, we synthesize an optimal policy for $\mathcal{M}$ by synthesizing an optimal policy for $\mathcal{M'}$.

 We use the nonlinear programming problem in \eqref{non_reach_objective}-\eqref{non_reach_cons6} to synthesize an optimal policy for $\mathcal{M'}$. 
\begin{subequations}
\label{unconstrained_program}
\begin{align}
\label{non_reach_objective}
& \underset{\lambda(s,a),\lambda(s)}{\text{maximize}}\ \   -\sum_{s\in S\backslash C}\sum_{t\in S}\eta(s,t)\log\Big(\frac{\eta(s,t)}{\nu(s)}\Big)\\
\label{non_reach_cons1}
& \text{subject to:}\nonumber\\
& \nu(s)-\sum_{t\in S\backslash C}\eta(t,s)=\alpha(s)  \ \quad  \forall s\in S\backslash C\\
\label{non_reach_cons2}
& \lambda(s) -\sum_{t\in S\backslash C}\eta(t,s)=\alpha(s) \   \quad \forall s\in C \\
\label{non_reach_cons3}
& \eta(s,t)=\sum_{a\in\mathcal{A}(s)}\lambda(s,a)\mathbb{P}_{s,a,t}  \quad   \forall t \in S, \forall s\in S\backslash C\\
\label{non_reach_cons4}
&\nu(s)=\sum_{a\in\mathcal{A}(s)}\lambda(s,a)  \qquad \qquad \forall s\in S\backslash C \\
\label{non_reach_cons5}
&\lambda(s,a)\geq 0 \qquad \qquad  \qquad   \qquad \forall a\in\mathcal{A}(s), \forall s\in S\backslash C\\
\label{non_reach_cons6}
&\lambda(s)\geq 0 \qquad \qquad \quad  \qquad  \qquad \forall s \in C
\end{align}
\end{subequations}

The decision variables in \eqref{non_reach_objective}-\eqref{non_reach_cons2} are $\lambda(s)$ for each $s$$\in$$C$, and $\lambda(s,a)$ for each $s$$\in$$S\backslash C$ and each $a$$\in$$\mathcal{A}(s)$. 
The function $\alpha$$:$$S$$\rightarrow$$\{0,1\}$ satisfies $\alpha(s_0)$$=$$1$ and $\alpha(s)$$=$$0$ for all $s$$\in$$S$$\backslash$$\{s_0\}$. Variables $\eta(s,t)$ and $\nu(s)$ are functions of $\lambda(s,a)$, and used just to simplify the notation. 

The constraints \eqref{non_reach_cons1}-\eqref{non_reach_cons2} represent the balance between the ``inflow" to and ``outflow" from states. The constraints \eqref{non_reach_cons3} and \eqref{non_reach_cons4} are used to simplify the notation and define the variables $\eta(s,t)$ and $\nu(s)$, respectively. The constraints \eqref{non_reach_cons5} and \eqref{non_reach_cons6} ensure that the expected residence time in the state-action pair $(s,a)$ and the probability of reaching the state $s$ is non-negative, respectively. We refer the reader to  \cite{Marta}, \cite{Puterman} for further details about the constraints.
 {\setlength{\parindent}{0cm}
\begin{prop} \label{prop_convex} The nonlinear program in \eqref{non_reach_objective}-\eqref{non_reach_cons6} is convex. $\quad \triangleleft$\end{prop}}

The above result indicates that a global maximum for the problem in \eqref{non_reach_objective}-\eqref{non_reach_cons6} can be computed efficiently. We now introduce Algorithm \ref{Algo_2} to synthesize an optimal policy for a given MDP with finite maximum entropy. 
{\setlength{\parindent}{0cm}
\begin{thm}
\label{Main_theorem_2}
 Let $\mathcal{M}$ be an MDP such that $H(\mathcal{M})$$<$$\infty$, $(C_i,D_i)$ $i$$=$$1,2,\ldots,n$ be MECs in $\mathcal{M}$, and $C$$:=$$\cup_{i=1}^nC_i$. For the input ($\mathcal{M}, C )$, Algorithm \ref{Algo_2} returns an optimal policy $\pi^{\star}$$\in$$\Pi^S(\mathcal{M})$ for $\mathcal{M}$, i.e. $H(\mathcal{M},{\pi^{\star}})$$=$$H(\mathcal{M})$. $\quad \triangleleft$
\end{thm}}

\begin{algorithm}[t]
 \caption{Synthesize the maximum entropy policy}
 \textbf{Require:} $\mathcal{M}$$=$$(S,s_0, \mathcal{A},\mathbb{P},\mathcal{AP},\mathcal{L})$ and $C$. \\
  \textbf{Return:} An optimal policy $\pi^{\star}$ for $\mathcal{M}$
  \begin{enumerate}[label=\arabic*:]
 \item Form the modified MDP $\mathcal{M'}$.
\item Solve \eqref{non_reach_objective}-\eqref{non_reach_cons6} for ($\mathcal{M'}$, $C$), and obtain $\lambda^{\star}(s,a)$.
\item \begin{algorithmic}
\For {$s$$\in$$ S$}
	\If {$s$$\not\in$$ C$}
		\If {$\sum_{a\in\mathcal{A}(s)}\lambda^{\star}(s,a)$$>$$0$}
		\State
			\State $\pi^{\star}_s(a)$$:=$$\frac{\lambda^{\star}(s,a)}{\sum_{a\in\mathcal{A}(s)}\lambda^{\star}(s,a)}$
		\Else
		 	\State  $\pi^{\star}_s(a)$$:=$$1$ for an arbitrary $a$$\in$$\mathcal{A}(s)$,
		\EndIf
	\Else
		\State $\pi^{\star}_s(a)$$:=$$1$ for an arbitrary $a$$\in$$\mathcal{A}(s)$.
	\EndIf
\EndFor
 \end{algorithmic}
 \end{enumerate}
 \label{Algo_2}
\end{algorithm} 
Proofs for above results can be found in Appendix \ref{proofs_appendix}. Computationally, the most expensive step of Algorithm \ref{Algo_2} is to solve the convex optimization problem \eqref{non_reach_objective}-\eqref{non_reach_cons6}. A solution whose objective value is arbitrarily close to the optimal value of \eqref{non_reach_objective} can be computed in time polynomial in the size of $\mathcal{M}$ via interior-point methods \cite{Serrano}, \cite{Nesterov}. Hence, the time complexity of Algorithm \ref{Algo_2} is polynomial in the size of $\mathcal{M}$.

\subsubsection{Unbounded maximum entropy} \label{unbounded_policy_synthesis} There is no optimal policy for this case due to \eqref{unbounded_definition_1}-\eqref{unbounded_definition_2}. Therefore, for a given MDP $\mathcal{M}$ and a constant $\ell$, we synthesize a policy $\pi$$\in$$\Pi^S(\mathcal{M})$ such that $H(\mathcal{M},\pi)$$\geq$$\ell$. Let $S_{B}$ be the union of all states in BSC MECs of $\mathcal{M}$, which can be found by using Algorithm \ref{algo_1}. We modify the MDP $\mathcal{M}$ by making all states $s$$\in$$S_B$ absorbing and denote the modified MDP by $\mathcal{M}'$. It can be shown that $H(\mathcal{M}')$$=$$H(\mathcal{M})$ by using arguments similar to the ones used in the proof of Proposition \ref{modified_MDP_prop}. As the first approach, we solve a convex feasibility problem. Specifically, we remove the objective in \eqref{non_reach_objective} and add the constraint
\begin{align}\label{residence_bound_L}
-\sum_{s\in S\backslash S_B}\sum_{t\in S}\eta(s,t)\log\Big(\frac{\eta(s,t)}{\nu(s)}\Big)\geq \ell
\end{align}
to the constraints in \eqref{non_reach_cons1}-\eqref{non_reach_cons6}. Then, we solve the resulting convex feasibility problem for ($\mathcal{M}'$, $S_B$, $\ell$) and obtain the desired policy $\pi$ by using the step 3 of Algorithm \ref{Algo_2}. 

Recall from Theorem \ref{Theorem1} that the unboundedness of the maximum entropy is caused by the existence of non-BSC MECs in $\mathcal{M}'$. In particular, we can induce MCs with arbitrarily large entropy by making the expected residence time in states contained in non-BSC MECs arbitrarily large. As the second approach, we bound the expected residence time in states $s$$\in$$S\backslash S_B$ in $\mathcal{M}'$ and relax this bound according to the desired level of entropy. Specifically, we add the constraint
\begin{align}\label{residence_bound}
\sum_{s\in S\backslash S_B}\sum_{a\in\mathcal{A}(s)}\lambda(s,a)\leq \Gamma
\end{align}
to the problem in \eqref{non_reach_objective}-\eqref{non_reach_cons6}. For the constraint \eqref{residence_bound}, $\Gamma$$\geq$$0$ is a predefined value and limits the expected residence time in states $s$$\in$$S\backslash S_B$. Let $H_{\Gamma}(\mathcal{M}')$ denote the maximum entropy $H(\mathcal{M}')$ of $\mathcal{M}'$ subject to the constraint \eqref{residence_bound}. Then, we have
\begin{align}\label{gamma_ineq}
H_{\Gamma}(\mathcal{M}')\geq H_{\Gamma'}(\mathcal{M}')
\end{align}
for $\Gamma$$\geq$$\Gamma'$, and $H_{\Gamma}(\mathcal{M}')$$=$$\infty$ for $\Gamma$$=$$\infty$. Therefore, by choosing an arbitrarily large $\Gamma$ value, we can synthesize a policy that induces an MC with arbitrarily large entropy.

\subsubsection{Infinite maximum entropy} \label{infinitee_policy_synthesis} 
The procedure to synthesize an optimal policy for MDPs with infinite maximum entropy is very similar to the ones described in Sections \ref{finite_policy_synthesis} and \ref{unbounded_policy_synthesis}. Therefore, we provide it in Appendix \ref{infinite_case_appendix}.
\section{Relating the maximum entropy of an MDP with the probability distribution of paths}\label{relate_paths}
In this section, we establish a link between the maximum entropy of an MDP $\mathcal{M}$ and the entropy of paths in an MC $\mathcal{M}^{\pi}$ induced from $\mathcal{M}$ by a stationary policy $\pi\in\Pi^S(\mathcal{M})$. 

\begin{figure}[b!]\vspace{-0.2cm}
\begin{subfigure}[b]{0.3\linewidth}
\scalebox{0.8}{
\begin{tikzpicture}[->, >=stealth', auto, semithick, node distance=2cm]

    \tikzstyle{every state}=[fill=white,draw=black,thick,text=black,scale=0.7]

    \node[state,initial,initial text=] (s_0) {$s_0$};
    \node[state] (s_1) [above right =6mm of s_0]  {$s_1$};
    \node[state] (s_2) [below right=6mm of s_0]  {$s_2$};
    \node[state] (s_3) [right=10mm of s_1]  {$s_3$};
    \node[state] (s_4) [right=10mm of s_2]  {$s_4$};

\path
(s_0)  edge     node{$a_1, 1$}     (s_1)
(s_0)	 edge      node{$a_2, 1$}     (s_2)
(s_1)	 edge       node{$a_1, 1$}     (s_3)
(s_1)	 edge       node{$a_2, 1$}     (s_4)
(s_2)	 edge       node{$ a_1,1$}     (s_4)
(s_3)	 edge    [loop right]       node{$a_1,1$}     (s_3)
(s_4)	 edge    [loop right]       node{$a_1,1$}     (s_4);
\end{tikzpicture}}
\caption{}
\label{fig:MDP_21}
\end{subfigure}
\hspace{0.18\linewidth}
\begin{subfigure}[b]{0.3\linewidth}
\scalebox{0.8}{
\begin{tikzpicture}[->, >=stealth', auto, semithick, node distance=2cm]

    \tikzstyle{every state}=[fill=white,draw=black,thick,text=black,scale=0.7]

    \node[state,initial,initial text=] (s_0) {$s_0$};
    \node[state] (s_1) [above right =6mm of s_0]  {$s_1$};
    \node[state] (s_2) [below right=6mm of s_0]  {$s_2$};
    \node[state] (s_3) [right=10mm of s_1]  {$s_3$};
    \node[state] (s_4) [right=10mm of s_2]  {$s_4$};

\path
(s_0)  edge     node{$2/3$}     (s_1)
(s_0)	 edge      node{$1/3$}     (s_2)
(s_1)	 edge       node{$1/2$}     (s_3)
(s_1)	 edge       node{$1/2$}     (s_4)
(s_2)	 edge       node{$1$}     (s_4)
(s_3)	 edge    [loop right]       node{$1$}     (s_3)
(s_4)	 edge    [loop right]       node{$1$}     (s_4);
\end{tikzpicture}}
\caption{}
\label{fig:MDP_22}
\end{subfigure}
\caption{(a) An MDP example \cite{Biondi}. (b) The MC induced by the policy that maximizes the entropy of the MDP. }
\label{fig:two-step-authenticate}
\end{figure}
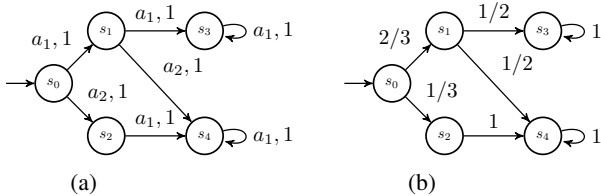

We begin with an example demonstrating the probability distribution of paths in an MC induced by a policy that maximizes the entropy of an MDP. Consider the MDP shown in Fig. \ref{fig:MDP_21} which is used in \cite{Biondi}. The policy that maximizes the entropy of the MDP is given by $\pi_{s_0}(a_1)$$=$$2/3$, $\pi_{s_0}(a_2)$$=$$1/3$, $\pi_{s_1}(a_1)$$=$$\pi_{s_1}(a_2)$$=$$1/2$. The MC induced by this policy is shown in Fig. \ref{fig:MDP_22}. There are three paths that reach the MECs, i.e., $(\{s_3\},\{a_1\})$ and $(\{s_4\},\{a_1\})$, of the MDP, each of which is followed with probability $1/3$ in the induced MC, i.e., the probability distribution of paths is uniform. 

Note that for the example given in Fig. \ref{fig:MDP_21}, the optimal policy that maximizes the entropy of the MDP is randomized, and action-selection at each state is performed in an online manner. In particular, an agent that follows the optimal policy chooses its action at each stage according to the outcomes of an online randomization mechanism. Therefore, it does not commit to follow a specific path at any state.

To rigorously establish the relation, illustrated in Fig. \ref{fig:MDP_21}, between the maximum entropy of an MDP and the entropy of paths in an induced MC, we need the following definitions.

 A \textit{strongly connected component} (SCC) $V$$\subseteq$$S$ in an MC $\mathcal{M}^{\pi}$ induced by a policy $\pi$$\in$$\Pi^S(\mathcal{M})$ is a maximal set of states in $\mathcal{M}^{\pi}$ such that for any $s$,$t$$\in$$V$, $(\mathcal{P}^{\pi})_{s,t}^n$$>$$0$ for some $n$$\in$$\mathbb{N}$. A \textit{bottom strongly connected component} (BSCC) $S_B$ in $\mathcal{M}^{\pi}$ is an SCC such that for all $s$$\in$$S_B$, $(\mathcal{P}^{\pi})_{s,t}^n$$=$$0$ for all $n$$\in$$\mathbb{N}$ and for all $t$$\in$$S\backslash S_B$.

In this section, for an induced MC $\mathcal{M}^{\pi}$, we denote the probability of a path with the finite path fragment $s_0\ldots s_n$ by
 \begin{align}\label{joint_prob_define}
 \mathcal{P}^{\pi}(s_0\ldots s_n):=\prod_{0\leq k < n} \mathcal{P}^{\pi}_{s_k, s_{k+1}},
 \end{align}
 which agrees with the probability measure introduced in Section \ref{Prelim}. Additionally, if the finite path fragment $s_0\ldots s_n$ in $\mathcal{M}^{\pi}$ satisfies $s_0,\ldots,s_{n-1}$$\not\in$$S_B$ and $s_n$$\in$$S_B$ for some $S_B$$\subseteq$$S$, we write $s_0\ldots s_n$$\in$$(S\backslash S_B)^{\star}S_B$. 
{\setlength{\parindent}{0cm}
\noindent \begin{definition}
(Entropy of paths) Let $\mathcal{M}^{\pi}$ be an MC induced from an MDP $\mathcal{M}$ by a stationary policy $\pi$$\in$$\Pi^S(\mathcal{M})$ and $S_B$$\subseteq$$S$ be the union of all BSCCs in $\mathcal{M}^{\pi}$. For $\mathcal{M}^{\pi}$, the \textit{entropy of the paths} that start from the initial state and reach a state in a BSCC in $\mathcal{M}^{\pi}$ is defined as
\begin{align}\label{entropy_paths_def}
&H(Paths^{\pi}(\mathcal{M})): =\nonumber \\
&\qquad \quad-\sum_{s_0\ldots s_n \in T } \mathcal{P}^{\pi}(s_0\ldots s_n)\log \mathcal{P}^{\pi}(s_0\ldots s_n)
\end{align}
where $T$$:=$$Paths_{fin}^{\pi}(\mathcal{M})$$\cap$$ (S\backslash S_B)^{\star}S_B$. 
\end{definition}}

A similar definition for the entropy of paths with fixed initial and final states can be found in \cite{Ekroot},\cite{Kafsi}. We note that
\begin{align}\label{sum_prob_paths}
\sum_{s_0\ldots s_n \in T } \mathcal{P}^{\pi}(s_0\ldots s_n)=1,
\end{align}
since any finite-state MC eventually reaches a BSCC \cite{Model_checking}. The following lemma establishes a relation between the entropy of paths and the entropy of an induced MC. 
{\setlength{\parindent}{0cm}
\begin{lemma}
\label{paths_lemma}
Let $\mathcal{M}$ be an MDP such that $H(\mathcal{M},{\pi})$$<$$\infty$ for any $\pi$$\in$$\Pi^S(\mathcal{M})$. Then, for any $\pi$$\in$$\Pi^S(\mathcal{M})$, we have
\begin{align}
H(Paths^{\pi}(\mathcal{M}))=H(\mathcal{M},{\pi}). \quad \triangleleft
\end{align}
\end{lemma}}\noindent

A proof for Lemma \ref{paths_lemma} can be found in Appendix \ref{proofs_appendix}. Finally, from the definition of the properties of the maximum entropy, Proposition \ref{memoryless_1} and Lemma \ref{paths_lemma}, we conclude that, if an MDP $\mathcal{M}$ has non-infinite maximum entropy, then we have
\begin{align}\label{equivalence_paths_entropy}
&H(\mathcal{M})\ = \ \sup_{\pi\in\Pi^S(\mathcal{M})}  \ H(Paths^{\pi}(\mathcal{M})).
\end{align}
The equality in \eqref{equivalence_paths_entropy} states that, for an MDP with non-infinite maximum entropy, a policy that maximizes the entropy of the MDP induces an MC with maximum entropy of paths among all MCs that can be induced from the MDP. Moreover, considering \eqref{sum_prob_paths}, such a policy maximizes the randomness of all paths with non-zero probability in an induced MC.
\section{Constrained Entropy Maximization for MDPs}\label{cons_section}
In this section, we consider the problem of maximizing the entropy of an MDP subject to an LTL constraint. We note that stationary policies are not sufficient to satisfy LTL constraints in general \cite{Model_checking}. Therefore, to be consistent with our definition of maximum entropy \eqref{max_ent_definition}, we first introduce the product MDP, over which LTL constraints are transformed into reachability constraints for which stationary policies are sufficient.  
\subsection{Product MDP}\label{product_section}
\par We construct an MDP that captures all paths of an MDP $\mathcal{M}$ satisfying an LTL specification $\varphi$ by taking the product of $\mathcal{M}$ and the DRA $A_{\varphi}$ corresponding to the specification $\varphi$.
{\setlength{\parindent}{0cm}
 \noindent \begin{definition} (Product MDP)
Let $\mathcal{M}$$=$$(S, s_0, \mathcal{A}, \mathbb{P}, \mathcal{AP}, \mathcal{L})$ be an MDP and $A_{\varphi}$$=$$(Q, q_0, 2^{\mathcal{AP}}, \delta, Acc)$ be a DRA. The product MDP $\mathcal{M}_p$$=$$\mathcal{M}$$\otimes$$ A_{\varphi}$$=$$(S_p, s_{0_p}, \mathcal{A}, \mathbb{P}_p, \mathcal{L}_p, Acc_p)$ is a tuple where
\begin{itemize}
\item $S_p$$=$$S $$\times$$ Q$,
\item $s_{0_p}=(s_0,q)$ such that $q=\delta(q_0,\mathcal{L}(s_0))$,
\item $\mathbb{P}_p((s,q), a, (s',q'))$=$\begin{cases} \mathbb{P}_{s,a,s'} & \text{if} \quad q'=\delta(q,\mathcal{L}(s')) \\ 0 & \text{otherwise}, \end{cases}$
\item $\mathcal{L}_p((s,q))=\{q\}$,
\item $Acc_p$$=$$\{(J_1^p,K_1^p),\ldots,(J_k^p,K_k^p) \}$ where $J_i^p$$=$$S$$\times$$J_i$ and $K_i^p$$=$$S$$\times$$K_i$ for all $(J_i,K_i)$$\in$$Acc$ and for all $i$$=$$1,\ldots, k$.
\end{itemize}
\end{definition}}
The product MDP $\mathcal{M}_p$ may contain unreachable states which can be found in time polynomial in the size of $\mathcal{M}_p$ by graph search algorithms, e.g., breadth-first search. Such states have no effect in the analysis of MDPs, and hence, can be removed from the MDP.  
We hereafter assume that there is no unreachable state in $\mathcal{M}_p$. 

There is a one-to-one correspondence between the paths of $\mathcal{M}_p$ and $\mathcal{M}$ \cite{Model_checking}. Moreover, a similar one-to-one correspondence exists between policies on $\mathcal{M}_p$ and $\mathcal{M}$. More precisely, for a given policy $\pi^{p}$$=$$\{\mu_0^p,\mu_1^p,\ldots\}$ on $\mathcal{M}_p$, we can construct a policy $\pi$$=$$\{\mu_0,\mu_1,\ldots\}$ on $\mathcal{M}$ by setting $\mu_i(s_i)$$=$$\mu_i^p((s_i,q_i))$. For a given policy $\pi^p$$\in$$\Pi^S(\mathcal{M}_p)$ on $\mathcal{M}_p$, the policy $\pi$$\in$$\Pi(\mathcal{M})$ constructed in this way is a non-stationary policy \cite{Model_checking}. 


Let $\pi^p$$\in$$\Pi^S(\mathcal{M}_p)$ be a policy on $\mathcal{M}_p$ and $\pi$$\in$$\Pi(\mathcal{M})$ be the policy on $\mathcal{M}$ constructed from $\pi^p$ through the procedure explained above. The paths of the MDP $\mathcal{M}$ under the policy $\pi$ satisfies the LTL specification $\varphi$ with probability of at least $\beta$, i.e., $\text{Pr}^{\pi}_{\mathcal{M}}(s_0$$\models$$\varphi)$$\geq$$\beta$, if and only if the paths of the product MDP $\mathcal{M}_p$ under the policy $\pi^p$ reaches accepting MECs in $\mathcal{M}_p$ with probability of at least $\beta$ and stays there forever \cite{Model_checking}. 
{\setlength{\parindent}{0cm}
\noindent \begin{definition} (Accepting MEC)
A MEC $(C,D)$ in a product MDP $\mathcal{M}_p$ with the set of accepting state pairs $Acc_p$ is an \textit{accepting} MEC if for some $(J^p,K^p)$$\in$$Acc_p$, $J^p$$\not\in$$C$ and $K^p$$\in$$C$.
\end{definition}}
Informally, accepting MECs are sets of states where the system can remain forever, and where the set $K_p$ is visited infinitely often and the set $J_p$ is visited finitely often. 

\vspace{-0.4cm}
\subsection{Constrained Problem}
 In this section, we formally state the constrained entropy maximization problem. Recall that, for an MDP $\mathcal{M}$, the problem of synthesizing a policy $\pi$$\in$$\Pi(\mathcal{M})$ that satisfies an LTL formula $\varphi$ with probability of at least $\beta$, i.e., $\text{Pr}^{\pi}_{\mathcal{M}}(s_0$$\models$$\varphi)$$\geq$$\beta$, is equivalent to the problem of synthesizing a policy $\pi$$\in$$\Pi^S(\mathcal{M}_p)$ that reaches the accepting MECs in $\mathcal{M}_p$ with probability of at least $\beta$ and stays there forever.

Our objective is to synthesize a policy that induces a stochastic process with maximum entropy whose paths satisfy the given LTL specification with desired probability. To this end, we synthesize a policy $\pi$$\in$$\Pi^S(\mathcal{M}_p)$ on $\mathcal{M}_p$ that reaches the accepting MECs in $\mathcal{M}_p$ with probability of at least $\beta$ and stays there forever.

We first partition the set $S_p$ of states of $\mathcal{M}_p$ into three disjoint sets as follows. We let $B$ be the set of all states in accepting MECs in $\mathcal{M}_p$, and $S_0$ be the set of all states that have zero probability of reaching the set $B$. Finally, we let $S_r$$=$$S_p\backslash \{B\cup S_0\}$ be the set of all states that are not in an accepting MEC in $\mathcal{M}_p$ and have nonzero probability of reaching the set $B$. These sets can be found in time polynomial in the size of $\mathcal{M}_p$ by graph search algorithms \cite{Model_checking}. 
 {\setlength{\parindent}{0cm}
\noindent \begin{problem}(\textbf{Constrained Entropy Maximization})\label{prob_2}
For a product MDP $\mathcal{M}_p$, verify whether there exists a policy $\pi^{\star}$$\in$$\Pi^S(\mathcal{M}_p)$ that solves the following problem:
\begin{subequations}
\begin{align}\label{max_ent_product_objective}
&\underset{\pi\in\Pi^S(\mathcal{M}_p)}{\text{maximize}}  \qquad H(\mathcal{M}_p,{\pi})\\ \label{LTL_constraint_product}
& \text{subject to:} \qquad  \text{Pr}^{\pi}_{\mathcal{M}_p}(s_0\models \lozenge B)\geq \beta
\end{align}
\end{subequations}
where $\text{Pr}^{\pi}_{\mathcal{M}_p}(s_0$$\models$$\lozenge B)$ denotes the probability of reaching the set $B$ from the initial state in $\mathcal{M}_p$ under the policy $\pi$. If such a policy exists, provide an algorithm to synthesize it. If it does not exist, provide a procedure to synthesize a policy $\pi'$$\in$$\Pi^S(\mathcal{M}_p)$ such that  $\text{Pr}^{\pi'}_{\mathcal{M}_p}(s_0\models \lozenge B)$$\geq$$\beta$ and $H(\mathcal{M},\pi')$$\geq$$\ell$ for a given constant $\ell$.
\end{problem}}

Note that if a policy that solves the problem in \eqref{max_ent_product_objective}-\eqref{LTL_constraint_product} chooses the actions in states $s$$\in$$B$ such that they form a BSCC in the induced MC, then the resulting policy ensures that the paths of the induced MC visit the states inside the set $B$ infinitely often and thus satisfies $\varphi$\cite{Model_checking}. 
\subsection{Policy synthesis}\label{product_policy_section}
  
In this section, for a product MDP $\mathcal{M}_p$ and its state partition $S_p$$=$$B$$\cup$$S_0$$\cup$$S_r$, we assume that $0$$<$$\beta$$\leq$$\max_{\pi\in\Pi(\mathcal{M})}\text{Pr}^{\pi}_{\mathcal{M}_p}(s_0$$\models$$\lozenge B)$, which can be verified in polynomial time by solving a linear optimization problem as shown in \cite{Model_checking, Marta}. We refer to a policy $\pi$$\in$$\Pi^S(\mathcal{M}_p)$ as an \textit{optimal policy} if it is a solution to the problem in \eqref{max_ent_product_objective}-\eqref{LTL_constraint_product} and chooses the actions in states $s$$\in$$B$ such that they form a BSCC in the induced MC.

For the synthesis of an optimal policy, we consider three cases according to the maximum entropy $H(\mathcal{M}_p)$ of $\mathcal{M}_p$, namely, finite, unbounded and infinite. 
 
 \subsubsection{Finite maximum entropy} \label{finite_constrained_case}
Let $(C_i,D_i)$ $i$$=$$1,2,\ldots,n$ be the MECs in $\mathcal{M}_p$, $C$$:=$$\cup_{i=1}^n C_i$, and $D$$:=$$\cup_{i=1}^n D_i$. We form the modified product MDP $\mathcal{M}_p'$ by making all states $s$$\in$$C$ absorbing in $\mathcal{M}_p$. We have $H(\mathcal{M}_p')$$=$$H(\mathcal{M}_p)$ due to Proposition \ref{modified_MDP_prop}. Recall that for a state $s$$\in$$C$, the variable $\lambda(s)$ in \eqref{non_reach_objective}-\eqref{non_reach_cons6} represents the probability of reaching the state $s$ from the initial state \cite{Marta}. Hence, we append the constraint
\begin{align}\label{reach_cons}
\sum_{s\in B}\lambda(s)\geq \beta
\end{align}
to the problem in \eqref{non_reach_objective}-\eqref{non_reach_cons6} in order to obtain a policy that induces an MC whose paths satisfy $\varphi$ with probability of at least $\beta$. Noting that $\beta$$\leq$$\sum_{s\in B}\lambda(s)$$\leq$$\max_{\pi\in\Pi(\mathcal{M})}\text{Pr}^{\pi}_{\mathcal{M}_p}(s_0$$\models$$\lozenge B)$, the resulting optimization problem always has a solution since its feasible set constitutes a closed compact set when the product MDP has finite maximum entropy.

The procedure to obtain a policy $\pi^{\star}_p$$\in$$\Pi^S(\mathcal{M}_p)$ that solves the problem in \eqref{max_ent_product_objective}-\eqref{LTL_constraint_product} for $\mathcal{M}_p$ with finite maximum entropy is as follows. First, we find MECs $(C,D)$ in $\mathcal{M}_p$ and form the modified MDP $\mathcal{M}_p'$ by making all states $s$$\in$$C$ absorbing. Second, we solve the problem in \eqref{non_reach_objective}-\eqref{non_reach_cons6} for $(\mathcal{M}_p', C, \beta)$ with the additional constraint \eqref{reach_cons}. Finally, we use step 3 of Algorithm \ref{Algo_2} to synthesize the policy $\pi^{\star}_p$$\in$$\Pi(\mathcal{M}_p)$. Note that the constructed policy ensures that, once reached, the system stays in the set $B$ forever, since all MECs in $\mathcal{M}_p$ with finite maximum entropy are bottom strongly connected. 

 \subsubsection{Unbounded maximum entropy}\label{unbounded_constrained_case} In this case, the product MDP $\mathcal{M}_p$ contains a non-BSC MEC due to Theorem \ref{Theorem1}. We assume that there is only one non-BSC MEC in $\mathcal{M}_p$, and it is contained in $S_r$. 
We first form the modified product MDP $\mathcal{M}'_p$ by making all states in BSC MECs in $\mathcal{M}_p$ absorbing. Note that $H(\mathcal{M}_p')$$=$$H(\mathcal{M}_p)$. Let $S_B$ denote the union of all absorbing states in $\mathcal{M}'_p$.
We verify the existence of a solution to the problem in \eqref{max_ent_product_objective}-\eqref{LTL_constraint_product} by solving the problem in \eqref{non_reach_objective}-\eqref{non_reach_cons6} for $(\mathcal{M}_p', S_B, \beta)$ with the additional constraint \eqref{reach_cons}. If the optimum value of the resulting problem is bounded, then we synthesize an optimal policy through step 3 of Algorithm \ref{Algo_2}. If it is not bounded, then there exists no optimal policy, in which case for a given constant $\ell$, we synthesize a policy $\pi^{\star}$$\in$$\Pi^S(\mathcal{M}_p)$ such that $H(\mathcal{M}_p,\pi^{\star})$$\geq$$\ell$ and $\text{Pr}^{\pi^{\star}}_{\mathcal{M}_p}(s_0\models \lozenge B)$$\geq$$ \beta$ by employing two different approaches.

As the first approach, we solve a convex feasibility problem. Specifically, for the problem in \eqref{non_reach_cons1}-\eqref{non_reach_cons6}, we remove the objective \eqref{non_reach_objective} and append the constraints \eqref{residence_bound_L} and \eqref{reach_cons} to the problem. Then, we solve the resulting convex feasibility problem for $(\mathcal{M}_p', S_B, \ell, \beta)$, and using step 3 of Algorithm \ref{Algo_2}, obtain a policy $\pi^{\star}_p$$\in$$\Pi^S(\mathcal{M}_p)$ such that $H(\mathcal{M}_p,\pi^{\star}_p)$$\geq$$\ell$ and $\text{Pr}^{\pi^{\star}_p}_{\mathcal{M}_p}(s_0$$\models$$\lozenge B)\geq \beta$. 

The second approach to obtain an induced MC with arbitrarily large entropy, whose paths satisfy the LTL specification with desired probability, is to bound the expected residence time in states $s$$\in$$S_p\backslash S_B$ and relax this bound according to the desired level of entropy. Specifically, we solve the problem in \eqref{non_reach_objective}-\eqref{non_reach_cons6} for $(\mathcal{M}_p', S_B, \beta, \Gamma)$ together with the constraints \eqref{residence_bound} and \eqref{reach_cons}, where $\Gamma$ is as defined in Section \ref{unbounded_policy_synthesis}. Then, by choosing an arbitrarily large $\Gamma$ value, we obtain an induced MC with the desired level of entropy. 

Finally, to ensure that the paths of the MC that is induced by the synthesized policy satisfies the LTL specification $\varphi$ with desired probability, we choose actions in states $s$$\in$$B$ such that $Succ(s)$$\subseteq$$B$.

 \subsubsection{Infinite maximum entropy}\label{infinitee_constrained_case}
For product MDPs with infinite maximum entropy, the verification of the existence and the synthesis of an optimal policy are achieved by procedures that are very similar to the ones presented in Sections \ref{finite_constrained_case} and \ref{unbounded_constrained_case}. Hence, we provide the analysis for product MDPs with infinite maximum entropy in Appendix \ref{infinite_case_appendix}. 

\section{Examples}\label{examples_section}
In this section, we illustrate the proposed methods on different motion planning scenarios. All computations are run on a 2.2 GHz dual core desktop with 8 GB RAM. All optimization problems are solved by using the splitting conic solver (SCS) \cite{SCS} in CVXPY \cite{cvxpy}.  For all LTL specifications, we construct deterministic Rabin automata using ltl2dstar \cite{ltl2dstar}.

In most motion planning scenarios, an agent can return to its current position by following different paths. Therefore, in general, the maximum entropy of an MDP that models the motion of an agent is either unbounded or infinite. However, as explained in Section \ref{policy_syntesis_section} and shown in the following examples, a policy that induces a stochastic process with an arbitrarily large entropy can easily be obtained by introducing constraints on the expected residence time in certain states. Additional motion planning examples are provided in \cite{Yagiz}. 
\subsection{Relation between entropy and exploration}
Randomizing an agent's paths while ensuring the completion of a task is important for achieving a better exploration of the environment \cite{Saerens} and obtaining a robust behavior against transition perturbations \cite{Deep_learning}. In this example, we demonstrate how the proposed method randomizes the agent's paths depending on the expected time until the completion of the task. 

\textit{Environment:} We consider the grid world shown in Fig. \ref{grid_graph} (left). The agent starts from the brown state. The red and green states are absorbing, i.e., once entered those states cannot be left. The agent has four actions in all other states, namely left, right, up and down. At each state, a transition to the chosen direction occurs with probability (w.p.) 0.7, and the agent slips to each adjacent state in the chosen direction w.p. 0.15. If the adjacent state in the chosen direction is a wall, e.g. up in brown state, a transition to the chosen direction occurs w.p. 0.85. If the state in the chosen direction is a wall, e.g., left in brown state, the agent stays in the same state w.p. 0.7 and moves to each adjacent state w.p. 0.15.  
\newcommand{\StaticObstacle}[2]{ \fill[red] (#1+0.1,#2+0.1) rectangle (#1+0.9,#2+0.9);}
\newcommand{\initialstate}[2]{ \fill[brown] (#1+0.15,#2+0.15) rectangle (#1+0.85,#2+0.85);}
\newcommand{\goalstate}[2]{ \fill[green] (#1+0.15,#2+0.15) rectangle (#1+0.85,#2+0.85);}
\begin{figure}[b]\vspace{-0.4cm}
\begin{subfigure}[b]{0.3\linewidth}
\scalebox{0.32}{
\begin{tikzpicture}
\draw[black,line width=1pt] (0,0) grid[step=1] (11,11);
\draw[black,line width=4pt] (0,0) rectangle (11,11);
			   \StaticObstacle{4}{3}  \StaticObstacle{5}{3}  \StaticObstacle{6}{3}
			   \StaticObstacle{4}{7}  \StaticObstacle{5}{7}  \StaticObstacle{6}{7}
			    \initialstate{0}{5}  \goalstate{10}{5}
			    
\node at (0.5,5+0.5) { \textbf{\huge S}};
\node at (10+0.5,5+0.5) { \textbf{\huge T}};
\node at (4+0.5,3+0.5) { \textbf{\huge B}};
\node at (5+0.5,3+0.5) { \textbf{\huge B}};
\node at (6+0.5,3+0.5) { \textbf{\huge B}};
\node at (4+0.5,7+0.5) { \textbf{\huge B}};
\node at (5+0.5,7+0.5) { \textbf{\huge B}};
\node at (6+0.5,7+0.5) { \textbf{\huge B}};

\end{tikzpicture}

}
\end{subfigure}
\hspace{0.2\linewidth}
\begin{subfigure}[b]{0.3\linewidth}
\scalebox{0.35}{
\begin{tikzpicture}
\draw[black,line width=1pt] (0,0) grid[step=1] (10,10);
\draw[black,line width=4pt] (0,0) rectangle (10,10);
			   \StaticObstacle{2}{2}  \StaticObstacle{3}{2}
			   \StaticObstacle{2}{3}  \StaticObstacle{3}{3}	
			 \StaticObstacle{6}{6}  \StaticObstacle{7}{6}
			 \StaticObstacle{6}{7}  \StaticObstacle{7}{7}  
			 \initialstate{0}{0} \goalstate{9}{9}
\node at (0.5,0.5) { \textbf{\huge S}};
\node at (5+0.5,0.5) { \textbf{\huge R1}};
\node at (8+0.5,3+0.5) { \textbf{\huge R2}};
\node at (1+0.5,6+0.5) { \textbf{\huge R3}};
\node at (4+0.5,8+0.5) { \textbf{\huge R4}};
\node at (9+0.5,9+0.5) { \textbf{\huge T}};
\node at (2+0.5,2+0.5) { \textbf{\huge B}};
\node at (2+0.5,3+0.5) { \textbf{\huge B}};
\node at (3+0.5,2+0.5) { \textbf{\huge B}};
\node at (3+0.5,3+0.5) { \textbf{\huge B}};
\node at (6+0.5,7+0.5) { \textbf{\huge B}};
\node at (7+0.5,7+0.5) { \textbf{\huge B}};
\node at (6+0.5,6+0.5) { \textbf{\huge B}};
\node at (7+0.5,6+0.5) { \textbf{\huge B}};
\end{tikzpicture}
}
\end{subfigure}
\caption{Grid world environments. The brown (S) and green (T) states are the initial and target states, respectively. The red (B) states are absorbing. }
\label{grid_graph}
\end{figure}
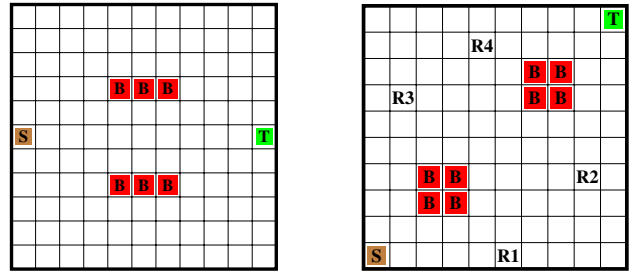

\textit{Task:} The agent's task is to reach and stay in the green state, labeled as $T$, while avoiding the red states, labeled as $B$. Formally, the task is $\varphi$$=$$\square \lnot B \land  \lozenge \square T$.

 We form the product MDP for the given task. It has 484 states, 1196 transitions, 10 MECs, and the average number of states in each MEC is 23. We require the agent to complete the task w.p. 1, i.e., $\text{Pr}_{\mathcal{M}}^{\pi}(s_0$$\models$$\varphi)$$=$$1$. The maximum entropy of the product MDP subject to the LTL constraint is unbounded. The minimum expected time $\Gamma$ required to complete the task $\varphi$ is roughly $14$ time steps, which can be computed by replacing the objective in \eqref{non_reach_objective} with ``minimize \ $\sum_{s\in S_r} \sum_{a\in \mathcal{A}(s)}\lambda(s,a)$" and appending \eqref{reach_cons} to the constraints in \eqref{non_reach_cons1}-\eqref{non_reach_cons6}. 

We synthesize two policies for two different expected times until the completion of the task. First, we synthesize a policy by requiring the agent to complete the task as fast as possible, i.e., $\Gamma$$=$$14$ time steps. Then, we synthesize a policy by allowing the agent to spend more time in the environment until the completion of the task, i.e., $\Gamma$$=$$60$ time steps. Solving the convex optimization problems take 122 and 166 seconds for $\Gamma$$=$$14$ and $\Gamma$$=$$60$ time steps, respectively.

\begin{figure}[t]
\begin{subfigure}[b]{0.3\linewidth}
\includegraphics[width= 3.95 cm, trim= {50 30 150 40 }, clip=true]{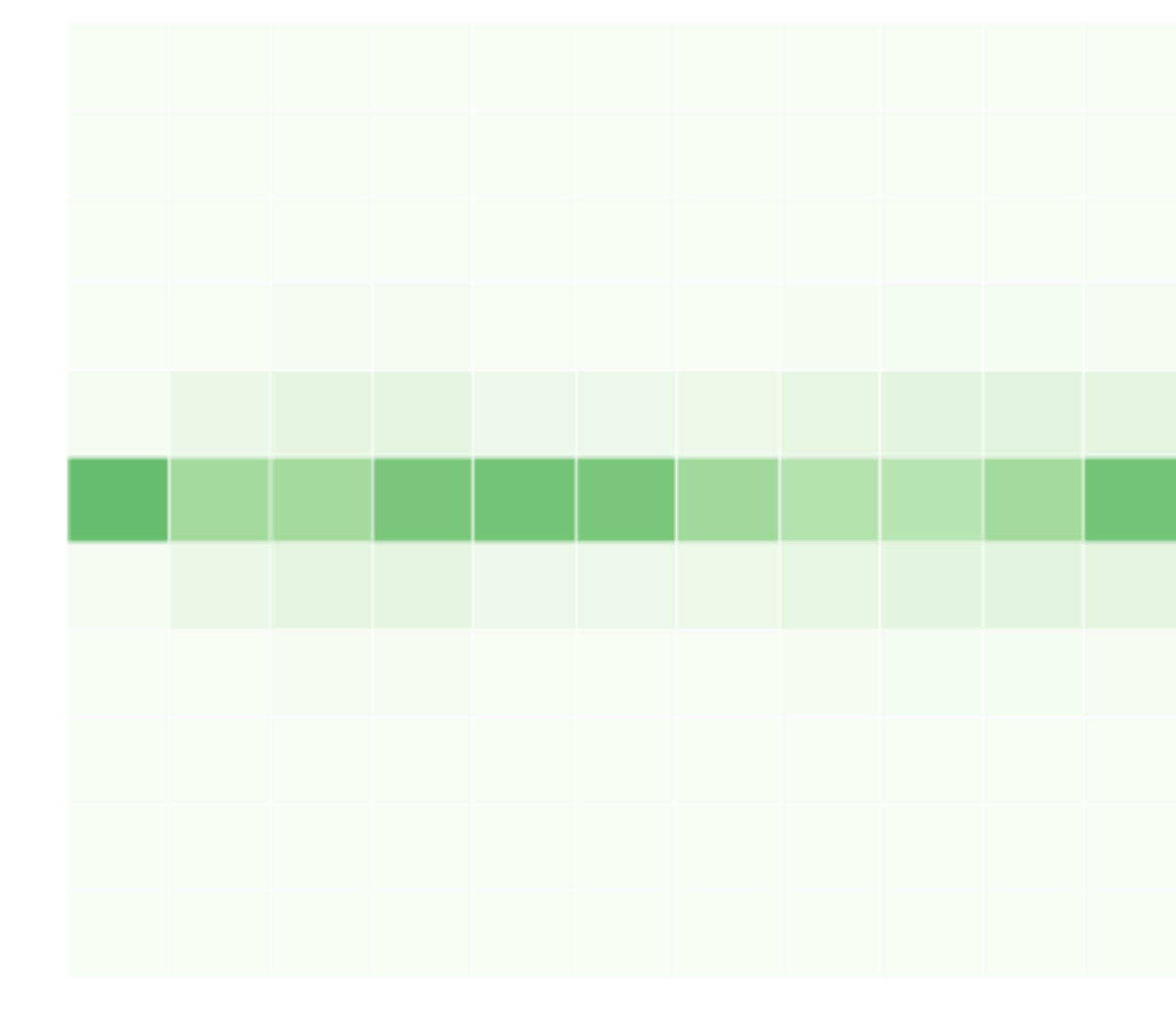}
\end{subfigure}
\hspace{0.15\linewidth}
\begin{subfigure}[b]{0.3\linewidth}
\includegraphics[width= 5 cm, trim= {50 30 50 40 }, clip=true]{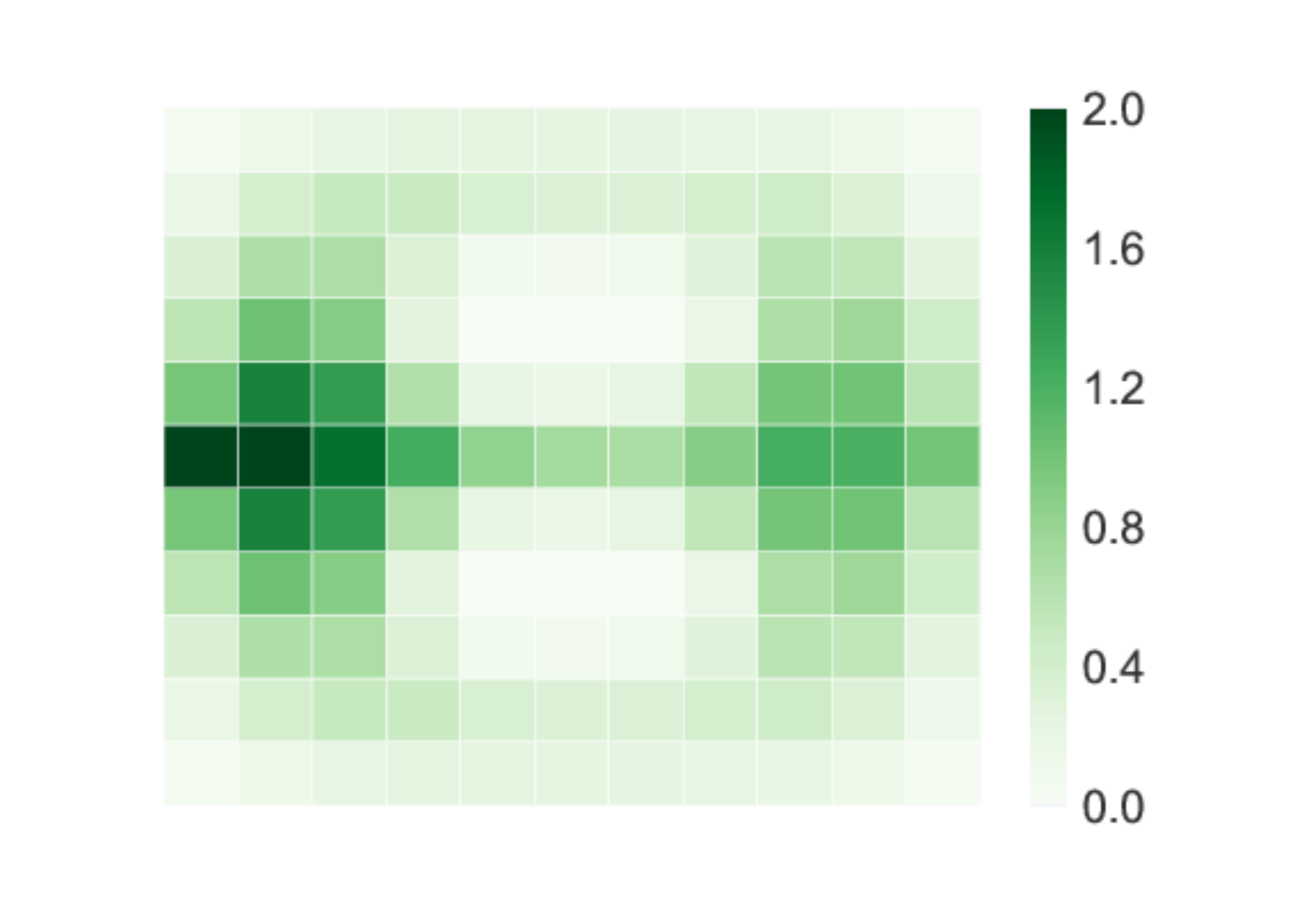}
\end{subfigure}
\caption{The expected residence time in states for different expected times $\Gamma$ until the completion of the task  \textit{in the same environment}. (Left) $\Gamma$$=$$14$ time steps, i.e., the minimum time required to complete the task with probability 1. (Right) $\Gamma$$=$$60$.   }
\vspace{-0.3cm}
\label{exploration_figure}
\end{figure}

The expected residence time in states for the induced MCs are shown in Fig. \ref{exploration_figure}.  We remind the reader that the environment is given in Fig. \ref{grid_graph} (left).
When the agent is given the minimum time $\Gamma$$=$$14$ time steps (left) to complete the task, it follows only the shortest paths, and therefore, cannot explore the environment. On the other hand, as it is allowed to spend more time, i.e., $\Gamma$$=$$60$ time steps (right), in the environment, it visits different states more often and utilizes different paths to complete the task. Consequently, the synthesized policy enables the continual exploration of the environment while ensuring the completion of the task.
\subsection{Relation between entropy and predictability}\label{example_1}

In this example, we consider an agent whose aim is to complete a task while leaking minimum information about its paths to an observer. We illustrate how the restrictions applied to the agent's paths by the task affect the predictability.

\textit{Environment:} We consider the grid world shown in Fig. \ref{grid_graph} (right). The agent starts from the brown (S) state. The red (B) states and green (T) state are absorbing. The agent has four actions in all other states, namely left, right, up and down. A transition to the chosen direction occurs w.p. 1 if the state in that direction is not a wall. If it is a wall, e.g., left direction in brown state, the agent stays in the same state w.p. 1.

\textit{Tasks:} We consider five increasingly restrictive task specifications for the agent which are listed in Table \ref{tasks_table}. The first task $\varphi_1$ is to reach and stay in the $T$ state while avoiding all red states. The second task $\varphi_2$ requires the agent to visit $R4$ state before completing the first task. The third task $\varphi_3$ requires the agent to visit $R3$ state before completing the second task and so on. 

\begin{table}[h!]\vspace{-0.3cm}
\centering
\caption{The agent's tasks. }
\scalebox{1.1}{
\begin{tabular}{ | l |}
  \hline			
$\varphi_1$$=$$\square \lnot Red \land \lozenge \square T$ \\  \hline  
 $\varphi_2$$=$$\square \lnot Red \land   \lozenge R4\land \lozenge \square T$ \\  \hline  
 $\varphi_3$$=$$\square \lnot Red \land \lozenge( R3\land \lozenge R4 )\land  \lozenge \square T$ \\  \hline  
 $\varphi_4$$=$$\square \lnot Red \land \lozenge( R2\land \lozenge (R3 \land \lozenge R4) )\land  \lozenge \square T$ \\  \hline  
 $\varphi_5$$=$$\square \lnot Red \land  \lozenge( R1\land \lozenge (R2 \land \lozenge (R3 \land \lozenge R4)) )\land   \lozenge \square T$ \\ \hline  
\end{tabular}}
\label{tasks_table}
\end{table}

\textit{Observer:} There is an observer that aims to predict the agent's paths in the environment. The observer is aware of the agent's task, knows the transition probabilities exactly, and runs yes-no probes in each state to determine the  successor state of the agent, i.e., probes that return an answer yes if the agent moves to the predicted successor state and no otherwise. The average number of yes-no probes run in a state is the expected number of observations needed by the observer to determine the correct successor state in that state \cite{Paruchuri}. The observer uses the Huffman procedure \cite{Huffman} to minimize the required number of probes. Let $\mathcal{P}_{s}$$=$$(\mathcal{P}_{s,1}, \mathcal{P}_{s,2}, \ldots, \mathcal{P}_{s,n} )$ be the transition probabilities from state $s$ to successor states sorted in decreasing order. The number of yes-no probes run in state $s$ is denoted by $\Upsilon_s$$=$$\mathcal{P}_{s,1}+\ldots+(n-1)\mathcal{P}_{s,n-1}$$+$$(n-1)\mathcal{P}_{s,n}$. The expected number of observations required to determine the agent's path is given by $O_{avg}$$=$$\sum_s\xi_s\Upsilon_s$, which weighs the required number of probes in each state with the expected residence time in the state. We refer the reader to \cite{Paruchuri} for further details about the observer model.

We form product MDPs for all tasks. The product MDP with the maximum number of states and transitions is the one constructed for the task $\varphi_5$. It has 800 states, 2138 transitions, 12 MECs, and the average number of states in each MEC is 29. For each task, we require the agent to complete the task w.p. 1. The maximum entropy of all product MDPs subject to corresponding LTL constraints are unbounded. We bound the expected time until the completion of any task by taking $\Gamma$$=$$33$ time steps, which is the minimum expected time required to complete the task $\varphi_5$, i.e., the most restrictive task. For each task, we synthesize a policy for the agent using the procedure explained in Section \ref{cons_section}. The longest computation time, which is for $\varphi_5$, is 15.2 seconds. 

\begin{figure}[h!]
\centering
\includegraphics[height=4 cm, width=7cm, trim= {0 0 0 30 }, clip=true]{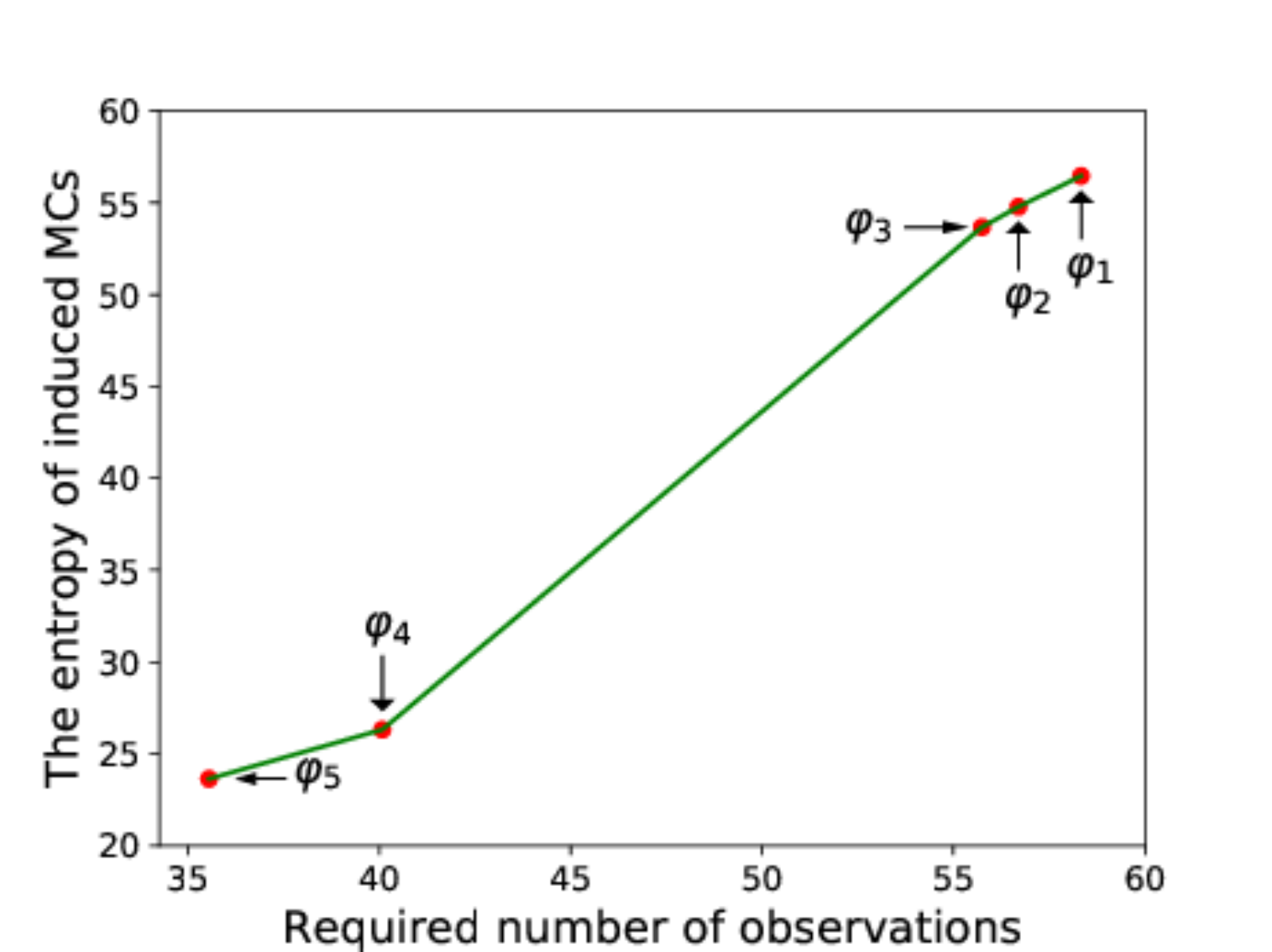}
\caption{ The relation between the maximum entropy of an MDP subject to an LTL constraint and the required number of observations to predict the agent's paths.    }
\label{tasks_figure}\vspace{-0.5cm}
\end{figure}

The entropy of Markov chains induced by the synthesized policies and the average number of observations required to predict the agent's paths are shown in Fig. \ref{tasks_figure}. As the task imposes more restrictions on the agent's paths, the entropy of the induced MC decreases and the prediction requires less observations in average. Additionally, as can be seen in Fig. \ref{tasks_figure}, the required numbers of observation for $\varphi_4$ and $\varphi_5$ are significantly smaller than those for $\varphi_1$, $\varphi_2$ and $\varphi_3$. This decrease is due to that these tasks force the agent to follow an ``S-shaped" path in a restricted time, i.e. $\Gamma=33$ time steps. For these tasks, although the agent still randomizes its paths to some extent, better predictability results cannot be obtained due to time restrictions.

\subsection{Predictability in a randomly generated MDP }
In this example, we investigate the relation between the probability of completing a task and the predictability of paths. We also evaluate the proposed algorithm against the algorithms introduced in \cite{Paruchuri}.

\textit{Environment:} We generate an MDP with 200 states, where each state has 8 randomly selected successor states. We choose four states, make them absorbing, and label three of them as ``unsafe" states and the remaining one as the ``target" state. The agent has 5 actions at each state, for which the transition probabilities to successor states are assigned randomly.  

\textit{Task:} The agent's task is to reach the target state while avoiding the unsafe states, i.e., $\varphi$$=$$\square \lnot unsafe \land \lozenge\square  target$.

\textit{Observer:} We use the same observer model introduced in Section \ref{example_1}.

\textit{Policies:} We compare the proposed method with weighted maximum entropy ($WME$) and binary search for randomization linear programming ($BRLP$) algorithms which are introduced in \cite{Paruchuri} for randomizing an agent's policy to minimize predictability. We note that in \cite{Paruchuri}, the authors claim that $WME$ algorithm is non-convex and cannot be solved in polynomial time. However, its convexity can be proven by Proposition \ref{prop_convex} since it solves a special case of the convex optimization problem given in \eqref{unconstrained_program}, i.e., it is equivalent to problem in \eqref{unconstrained_program} when transition probabilities are either 0 or 1. We refer the reader to \cite{Paruchuri} for further details about the $WME$ and $BRLP$ algorithms. 

 We form the product MDP. It has 800 states, 2172 transitions and 5 MECs each of which contains a single state. The maximum probability of completing the task $\varphi$ is obtained as $\beta$$=$$0.75$ by solving a linear programming problem introduced in \cite{Marta}. The maximum entropy of the product MDP subject to the LTL constraint $\text{Pr}_{\mathcal{M}}^{\pi}(s_0$$\models$$\varphi)$$\geq$$\beta$ is unbounded for all $\beta$$>$$0$. We fix the expected time until the completion of the task to $\Gamma$$=$$200$ time steps, and synthesize policies for different values of $\beta$. Solving the optimization problems takes at most 150, 155, and 92 seconds for the proposed method, $WME$ and $BRLP$ algorithms, respectively. 

\begin{figure}[h!]
\centering
\includegraphics[height=4 cm, width=7cm, trim= {0 0 0 30 }, clip=true]{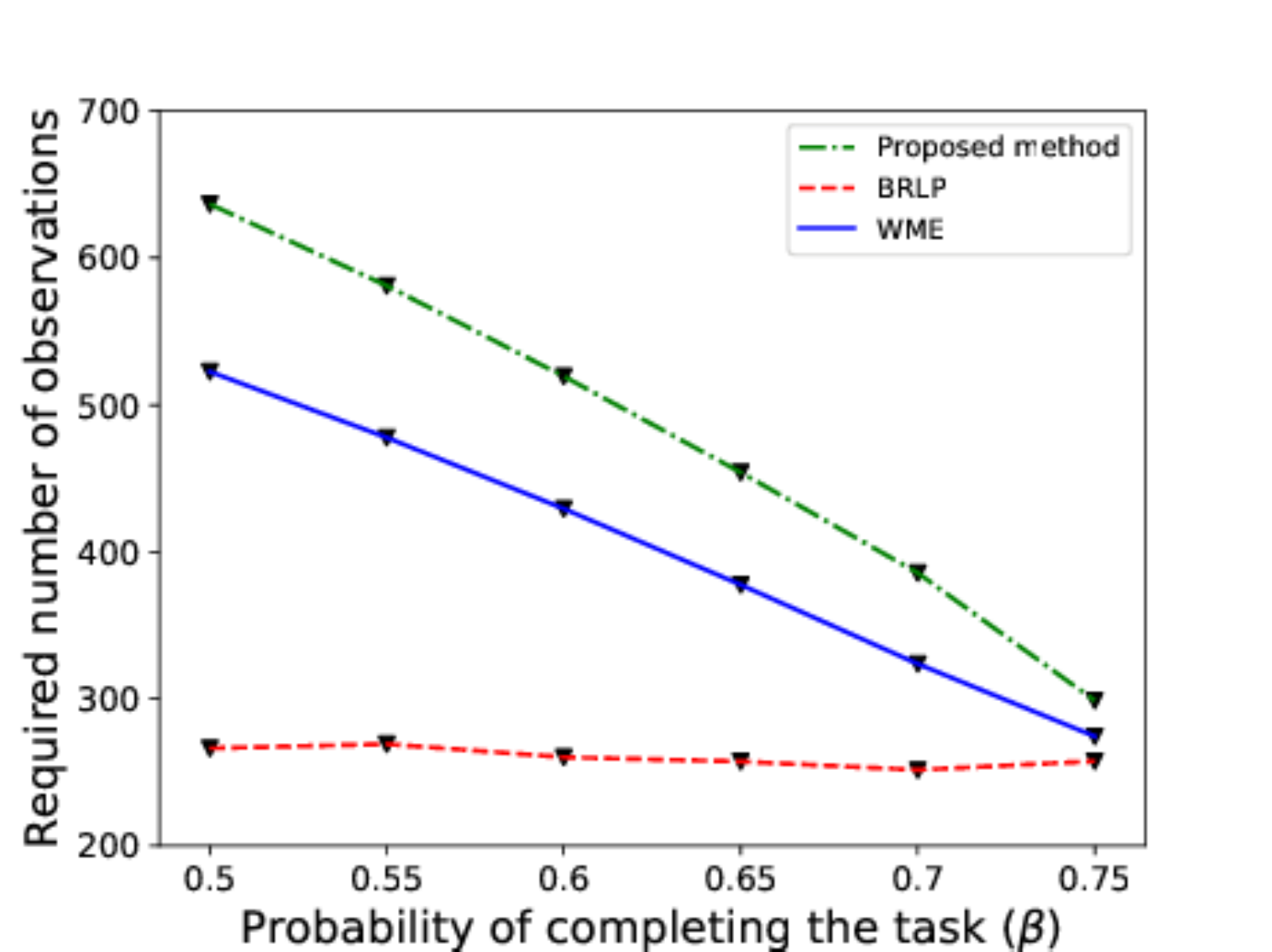}
\caption{The trade-off between the probability of completing the task and the predictability of paths. $BRLP$ and $WME$ are algorithms proposed in \cite{Paruchuri} to randomize the agent's actions.}
\label{compare_paruchuri}
\end{figure}

The required number of observations to predict the agent's paths for different $\beta$ values are shown in Fig. \ref{compare_paruchuri}. As the probability of completing the task decreases, the randomness of the agent's paths increases and the prediction requires more observations in average. Therefore, there is a trade-off between the probability of satisfying the task and the randomness of the paths. Additionally, the proposed method (green) requires two times more observations than the $BRLP$ algorithm (red) when $\beta$$=$$0.5$. Note also that the $WME$ algorithm (blue) cannot achieve better predictability results than the proposed method because it does not exploit the inherent stochasticity in the environment and rely solely on the randomization of the agent's actions to generate unpredictable paths.

\section{Conclusions and Future Work}\label{conclusion_section}
We showed that the maximum entropy of an MDP can be either finite, infinite or unbounded, and presented an algorithm to verify the property of the maximum entropy for a given MDP. We presented an algorithm, based on a convex optimization problem, to synthesize a policy that maximizes the entropy of an MDP. For MDPs with non-infinite maximum entropy, we established the equivalence between the maximum entropy of an MDP and the maximum entropy of paths in the MDP. Finally, we provided a procedure to obtain a policy that maximizes the entropy of an MDP while ensuring the satisfaction of a temporal logic specification with desired probability.

An interesting future direction is to include adversaries to the framework by modeling the problem as a two-player game. Being informed about the aims and capabilities of rational/irrational adversaries in the environment, an agent may want to explore its environment while avoiding the threats caused by adversaries. Another future direction may be to extend this work to multi-agent scenarios by describing the tasks, and communication and coordination constraints between the agents as temporal logic specifications.

\bibliographystyle{IEEEtran}
\bibliography{FinalVersion}


\appendices

\section{}\label{proofs_appendix}
We first define the hitting probability of a set of states in an MC $\mathcal{M}^{\pi}$ induced from an MDP $\mathcal{M}$ by a policy $\pi$$\in$$\Pi^S(\mathcal{M})$. If a sequence $s_k s_{k+1}\ldots s_n$ in $\mathcal{M}^{\pi}$ of states satisfies $s_k$$\in$$B$, $s_{k+1},\ldots s_{n-1}$$\not\in$$A$ and $s_n$$\in$$A$, we write  $s_k s_{k+1}\ldots s_n$$\in$$B (S\backslash A)^{\star} A$.
{\setlength{\parindent}{0cm}
\noindent \begin{definition}
 For an induced MC $\mathcal{M}^{\pi}$, the \textit{hitting probability} of a set $A$ of states from a set $B$ of states is $\rho^{\pi}(B,A)$$:=$$\sum_{s_k\ldots s_n\in T}\prod_{k\leq i < n} \mathcal{P}^{\pi}_{s_i, s_{i+1}}$ where $T$$=$$B (S\backslash A)^{\star} A$.
\end{definition}}

\noindent \textbf{Proof of Proposition \ref{memoryless_1}} If $\sup_{\pi\in\Pi^S(\mathcal{M})}H(\mathcal{M},{\pi})$$=$$\infty$, equality follows from the fact that $\Pi^S(\mathcal{M})$$\subseteq$$\Pi(\mathcal{M})$. If $\sup_{\pi\in\Pi^S(\mathcal{M})}H(\mathcal{M},{\pi})$$<$$\infty$, the result follows from Proposition \ref{modified_MDP_prop} in this paper, Proposition 35 in \cite{Chen}, and Proposition 2 in \cite{Bertsekas_SSP}. Specifically, $H(\mathcal{M},\pi)$ can be written as an expected total cost with respect to a specific cost function on an MDP with a compact action set \cite{Chen}, Proposition \ref{modified_MDP_prop}, together with Theorem \ref{Theorem1}, implies that every stationary policy on this MDP is \textit{proper}, i.e., all stochastic processes induced by stationary policies are guaranteed to reach an absorbing state within finite time step, and finally, the sufficiency of stationary policies to minimize the expected total cost on this MDP follows from \cite{Bertsekas_SSP}.$\Box$

\noindent\textbf{Proof of Lemma \ref{Successor}:} 
We prove the sufficiency by contradiction and the necessity by construction.

($\Rightarrow$) To obtain a contradiction, assume there exists a state $s$$\in$$C$ that satisfies $\lvert \cup_{a\in D(s)} Succ(s,a)\rvert$$=$1 and is both stochastic and recurrent in an MC $\mathcal{M}^{\pi}$ induced by a stationary policy $\pi$$\in$$\Pi^S(\mathcal{M})$. Since the state $s$ is stochastic in $\mathcal{M}^{\pi}$ and $\lvert \cup_{a\in D(s)} Succ(s,a)\rvert$$=$$1$, the policy $\pi$ satisfies $\pi_s(a_2)$$>$0 for some action $a_2$$\in$$D_0(s)$. Therefore, there exists a state $t$ in $\mathcal{M}^{\pi}$ such that $t$$\in$$Succ(s)$$\backslash$$C$ by the definition of $D_0(s)$.

Let $u$$=$$\mathcal{P}^{\pi}_{s,t}$$>$0. Then, for all $k$$\in$$\mathbb{N}$, $(\mathcal{P}^{\pi})_{s,s}^k$$\leq($1$-$$u')^k$ for some $0$$<$$u'$$\leq$$u$. (Note that if $u'$$=$$0$ for some $k$$\in$$\mathbb{N}$, then there exists a path that starts from the state $s$, visits the state $t$ and returns to the state $s$ with probability 1. However, in this case we should have $t$$\in$$C$.) As a result, 
\begin{align}
\label{proof_lemma_2}
\xi^{\pi}_s=\rho^{\pi}(s_0,s)\sum_{k=0}^{\infty}(\mathcal{P}^{\pi})_{s,s}^k\leq\sum_{k=0}^{\infty}(1-u')^k=\frac{1}{u'}<\infty,
\end{align}
where we use the fact that the hitting probability satisfies $\rho^{\pi}(s_0,s)$$\leq$$1$. This raises a contradiction since the state $s$ is recurrent, and it must satisfy $\xi^{\pi}_s$$=$$\infty$.

($\Leftarrow$) Suppose there exists a state $s$$\in$$ C$ such that $\lvert \cup_{a\in D(s)} Succ(s,a)\rvert$$>$1. Then, either (i) there exist actions $a_i$$,$$a_j$$\in$$ D(s)$ such that $Succ(s,a_i)$$\backslash$$ Succ(s,a_j)$$\neq$$\emptyset$, or (ii) there exists an action $a_i$$\in$$ D(s)$ such that $|Succ(s,a_i)|$$>$$1$. For case (i), we construct a policy $\pi$$\in$$\Pi^S(\mathcal{M})$ such that $\pi_s(a_i)$$>$0 and $\pi_s(a_j)$$>$$0$ in the state $s$, and $\pi_{s'}(a)$$=$$1$ for some action $a$$\in$$D(s')$ in states $s'$$\in$$C\backslash \{s\}$. Note that such actions exist by the definition of MEC. Finally, in states $t$$\not\in$$C$, we choose actions such that the state $s$ is reachable from the initial state. In the MC $\mathcal{M}^{\pi}$ induced by $\pi$, the state $s$ is an element of a bottom strongly connected component (BSCC) and $|Succ(s)|$$>$$1$. Hence, it is both recurrent and stochastic. For the case (ii), we choose a policy $\pi$$\in$$\Pi^S(\mathcal{M})$ such that $\pi_s(a_i)$$=$1 in the state $s$ and $\pi_{s'}(a)$$=$$1$ for some $a$$\in$$D(s')$ in states $s'$$\in$$C\backslash \{s\}$. In $\mathcal{M}^{\pi}$, the state $s$ belongs to a BSCC and has multiple successor states. Hence, it is both recurrent and stochastic. $\Box$\\
\noindent \textbf{Proof of Theorem \ref{Theorem1}}:
We first prove the necessary and sufficient conditions for an MDP to have infinite or unbounded maximum entropy. Then, we show that if the maximum entropy is not infinite and not unbounded, then it is finite and attainable by a stationary policy.

\noindent \textbf{Infinite maximum entropy.} We prove that the maximum entropy of an MDP $\mathcal{M}$ is infinite if and only if there exists a state $s$$\in$$C$ such that $\lvert \cup_{a\in D(s)} Succ(s,a)\rvert$$>$1, and conclude, by Lemma \ref{Successor}, that the claim holds.

 ($\Rightarrow$) The proof is by contradiction. Assume that the maximum entropy of $\mathcal{M}$ is infinite, i.e. $\max_{\pi\in\Pi^S(\mathcal{M})}H(\mathcal{M},{\pi})$$=$$\infty$, and $\lvert \cup_{a\in D(s)}Succ(s,a)\rvert$$=$1 for all states $s$$\in$$C$. We consider two cases: (i) $D_0(s)$$=$$\emptyset$ for all $s$$\in$$C$, and  (ii) $D_0(s)$$\neq$$\emptyset$ for some $s$$\in$$C$.

\textbf{Case (i)}: Suppose that $D_0(s)$$=$$\emptyset$ for all $s$$\in$$C$. Then, for an arbitrarily chosen MC $\mathcal{M}^{\pi}$ induced by a policy $\pi$$\in$$\Pi^S(\mathcal{M})$, we have $\lvert Succ(s)\rvert $$=$$1$ for all $s$$\in$$C$. Hence, for $\mathcal{M}^{\pi}$, $L^{\pi}(s)$$=$$0$ for all $s$$\in$$C$ due to \eqref{local_entropy_def}. Recall that $H(\mathcal{M},{\pi})$$<$$\infty$ if and only if for all $s$$\in$$S$, $\xi^{\pi}_s$$=$$\infty$ implies $L^{\pi}(s)$$=$$0$ \cite{Biondi}, and note that, if $\xi^{\pi}_s$$=$$\infty$, then $s$$\in$$C$. Consequently, we have $H(\mathcal{M},{\pi})$$<$$\infty$ for any MC $\mathcal{M}^{\pi'}$ induced by a policy $\pi'$$\in$$\Pi^S(\mathcal{M})$ since we choose $\mathcal{M}^{\pi}$ arbitrarily. This implies that $\max_{\pi\in\Pi^S(\mathcal{M})}H(\mathcal{M},{\pi})$$<$$\infty$ and raises a contradiction since for an MDP with infinite maximum entropy, we have $\max_{\pi\in\Pi^S(\mathcal{M})}H(\mathcal{M},{\pi})$$=$$\infty$.

\textbf{Case (ii):} Suppose that $D_0(s)$$\neq$$\emptyset$ for some $s$$\in$$C$. For an arbitrarily chosen induced MC $\mathcal{M}^{\pi}$, the local entropy of any state $s$$\in$$S$ is bounded by $L^{\pi}(s)$$\leq$$\log\lvert S\rvert$ \cite{Biondi}. We assume, without loss of generality, that for all states $s$$\in$$C$, $\lvert Succ(s)\rvert $$>$$1$ for $\mathcal{M}^{\pi}$. (If $\lvert Succ(s)\rvert $$=$$1$, the state $s$ has no contribution to the entropy of $\mathcal{M}^{\pi}$ due to \eqref{process_entropy}.) Recalling \eqref{proof_lemma_2}, for any $s$$\in$$C$, there exists a constant $u_s$$>$$0$ such that $\xi^{\pi}_s$$\leq$$\frac{1}{u_s}$. Then, for any $\mathcal{M}^{\pi}$, $\xi^{\pi}_s$$\leq$$\frac{1}{u'}$ for all $s$$\in$$C$, where $u'$$=$$\min \{u_s$$:$$\ \xi^{\pi}_s$$\leq$$\frac{1}{u_s},\ s$$\in$$ C \}$.

We now consider the states $s'$$\not\in$$C$. For any state $s'$$\not\in$$C$,  $\mathbb{P}_{s'as'}$$<$$1$ for all $a$$\in$$\mathcal{A}(s')$ since otherwise the state $s'$ must belong to a MEC. Then, there exists a constant $u''$$>$$0$ such that $\xi^{\pi}_{s'}$$\leq$$\frac{1}{u''}$ for all $s'$$\not\in$$C$ for any induced MC $\mathcal{M}^{\pi}$. As a result, 
\begin{align*}
H(\mathcal{M},{\pi})&=\sum_{s\in S}L^{\pi}(s)\xi^{\pi}_s=\sum_{s\in C}L^{\pi}(s)\xi^{\pi}_s+\sum_{s\not\in C}L^{\pi}(s)\xi^{\pi}_s,\\
&\leq\frac{\lvert S\rvert \log\lvert S\rvert}{u'}+\frac{\lvert S\rvert \log\lvert S\rvert}{u''}< \infty,
\end{align*}
\noindent for any induced MC $\mathcal{M}^{\pi}$ and for some $u',u''$$>$$0$. This implies $\max_{\pi\in\Pi^S(\mathcal{M})}H(\mathcal{M},{\pi})$$<$$\infty$ and raises a contradiction.

($\Leftarrow$) Using the proof of Lemma \ref{Successor}, we can construct a policy for which the induced MC contains a state that is both stochastic and recurrent. By Corollary 1 in \cite{Biondi}, the entropy of the induced MC is infinite. \\
\textbf{Unbounded maximum entropy.} ($\Rightarrow$) The proof is by contradiction. Assume that the maximum entropy of $\mathcal{M}$ is unbounded, and there exists $s$$\in$$C$ such that $\lvert \cup_{a\in D(s)}Succ(s,a)\rvert$$>$$1$ \textit{or} all MECs in $\mathcal{M}$ are bottom strongly connected. First, suppose that $H(\mathcal{M})$ is unbounded and there exists $s$$\in$$C$ such that $\lvert \cup_{a\in D(s)}Succ(s,a)\rvert$$>$$1$. Then, by case (i) of Theorem \ref{Theorem1}, the maximum entropy of $\mathcal{M}$ is infinite, which is a contradiction. Second, suppose that $H(\mathcal{M})$ is unbounded and all MECs in $\mathcal{M}$ are bottom strongly connected. Then, $H(\mathcal{M},{\pi})$$<$$\infty$ for all $\pi$$\in$$\Pi^S(\mathcal{M})$ by the definition of unboundedness. Using case (i) of Theorem \ref{Theorem1}, we conclude that there is no state in MECs that is both stochastic and recurrent in an induced MC $\mathcal{M}^{\pi}$. Consequently, all states $s$$\in$$C$ are deterministic for any induced MC since $D_0(s)$$=$$\emptyset$ and $\lvert \cup_{a\in D(s)}Succ(s,a)\rvert$$=$$1$ for all $s$$\in$$C$. This implies that $L^{\pi}(s)$$=$$0$ for all $s$$\in$$C$. Since every state $s'$$\not\in$$C$ satisfies $\mathbb{P}_{s',a,s'}$$<$$1$ for all $a$$\in$$\mathcal{A}(s')$, there exists a constant $u^{\star}$$>$0 such that for all $s'$$\not\in$$C$ and for all $\pi$$\in$$\Pi^S(\mathcal{M})$, $\xi^{\pi}_{s'}$$\leq$$\frac{1}{u^{\star}}$. As a result,
\begin{align*}
H(\mathcal{M},{\pi})=\sum_{s\in S}L^{\pi}(s)\xi^{\pi}_s=\sum_{s\not\in C}L^{\pi}(s)\xi^{\pi}_s\leq\frac{\lvert S\rvert \log\lvert S\rvert}{u^{\star}}
\end{align*}
for any policy $\pi$$\in$$\Pi^S(\mathcal{M})$. Hence, the maximum entropy is bounded. Since we assumed at the beginning that the maximum entropy is unbounded, this raises a contradiction.

($\Leftarrow$) The proof is by construction. Suppose that $\mathcal{M}$ has a MEC $C^{\star}$ which is not bottom strongly connected. Then, there exists a state $s$$\in$$ C^{\star}$ such that $D_0(s)$$\neq$$\emptyset$. Let $R$$=$$Succ(s,a_i)\backslash C^{\star}$ for $a_i$$\in$$D_0(s)$, i.e., the set of states that are reachable from the state $s$ and do not belong to the MEC $C^{\star}$. We construct a policy $\pi$$\in$$\Pi^S(\mathcal{M})$ such that, for the state $s$, $\sum_{t\in R}\mathcal{P}^{\pi}_{s,t}$$=$$\epsilon$, and for states $s'$$\in$$C^{\star}\backslash \{s\}$, $\pi_{s'}(a)$$=$$1$, for some $a$$\in$$D(s')$. For states $s{''}\not\in C^{\star}$, we choose actions so that state $s$ is reachable from the initial state in the induced MC $\mathcal{M}^{\pi}$. 

The induced MC $\mathcal{M}^{\pi}$ has the property that $\rho^{\pi}(t,s)$$=$$1$ for all $t$$\in$$C^{\star}\backslash \{s\}$ and $\rho^{\pi}(s,s)$$=$$1$$-$$\epsilon'$ for some $0$$<$$\epsilon'$$\leq$$\epsilon$. Here, we note that state $s$ is not recurrent in $\mathcal{M}^{\pi}$.

Since states $s$$\in$$C^{\star}$ are reachable from the initial state in $\mathcal{M}^{\pi}$, we have $\rho^{\pi}(s_0,s)$$>$$0$. Additionally, $\mathcal{K}$$:=$$\{k | (\mathcal{P}^{\pi})_{s_0,s}^k$$>$$0 \}$ is non-empty. Let $k^{\star}$$:=$$\min(\mathcal{K})$ and $\rho$$:=$$(\mathcal{P}^{\pi})_{s_0,s}^{k^{\star}}$. Then, we have $\rho$$\leq$$\rho(s_0,s)$ because $\rho$ only includes the first hitting probability. Moreover, the state $s$ satisfies  $\rho^{\pi}(s,s)$$=$$1$$-$$\epsilon'$$<$$1$. Then, $\xi_s^{\pi}$$=$$\frac{\rho^{\pi}(s_0,s)}{1-\rho^{\pi}(s,s)}$$\geq$$\frac{\rho}{\epsilon'}$$\geq$$\frac{\rho}{\epsilon}$, where the equality is a well-known result for finite-state MCs \cite{Residence_time}, \cite{Serfozo}.

The local entropy of the state $s$ is the smallest when it has two outgoing transitions, one with probability $\epsilon$ to a state $t$$\in$$R$ and the other with probability $1$$-$$\epsilon$ to a state in $C^{\star}$ \cite{Biondi}. (It can be imagined as a Bernoulli random variable with parameter $\epsilon$ where $\epsilon$ can be arbitarily small.). Hence, $L^{\pi}(s)$$\geq$$-$$((\epsilon\log\epsilon)$$+$$((1$$-$$\epsilon)\log(1$$-$$\epsilon)))$. As a result, 
\begin{align}
L^{\pi}(s)\xi^{\pi}_s\geq -\frac{\rho((\epsilon\log\epsilon)+((1-\epsilon)\log(1-\epsilon)))}{\epsilon}.
\end{align}
Note that $\lim_{\epsilon\rightarrow 0^+}L^{\pi}(s)\xi^{\pi}_s$$=$$\infty$. Therefore, for any policy $\pi$$\in$$\Pi^S(\mathcal{M})$, it is always possible to find another policy $\pi'$$\in$$\Pi^S(\mathcal{M})$ that induces an MC with a greater entropy. Hence, the maximum entropy of the MDP is unbounded. 

\noindent\textbf{Finite maximum entropy.} ($\Leftarrow$) The result follows from the definition of the finite maximum entropy.

($\Rightarrow$) Assume that the maximum entropy is not infinite and not unbounded. Hence, $H(\mathcal{M})$$=$$\sup_{\pi\in\Pi(\mathcal{M})}H(\mathcal{M},\pi)$$<$$\infty$ which implies that, for any policy $\pi$$\in$$\Pi^S(\mathcal{M})$, $H(\mathcal{M},{\pi})$$<$$\infty$. Then, for all states $s$$\in$$C$, $\lvert \cup_{a\in D(s)}Succ(s,a)\rvert$$=$$1$ and all MECs are BSC by cases (i) and (ii) of Theorem \ref{Theorem1}, respectively. As a result, for all $s$$\in$$C$, we have $L^{\pi}(s)$$=$$0$ for any $\pi$$\in$$\Pi^S(\mathcal{M})$, and hence, $H(\mathcal{M},{\pi})=\sum_{s\in S\backslash C}\xi^{\pi}_sL^{\pi}(s)$ for any $\pi$$\in$$\Pi^S(\mathcal{M})$.

Suppose that there exists a state $s$$\in$$S\backslash C$ such that $\xi^{\pi}_s$$=$$0$ for some $\pi$$\in$$\Pi^S(\mathcal{M})$. Then, for the induced MC $\mathcal{M}^{\pi}$, $\xi^{\pi}_sL^{\pi}(s)$$=$$0$ since $L^{\pi}(s)$$\leq$$\log\lvert S \rvert$ is bounded. Therefore, without loss of generality, we can neglect unreachable states in any induced MC $\mathcal{M}^{\pi}$ and assume $\xi^{\pi}_s$$>$$0$ for all states $s$$\in$$S\backslash C$. We pick an arbitrary state $s$$\in$$S\backslash C$ and an arbitrary policy $\pi$$\in$$\Pi^S(\mathcal{M})$, and define a new function $\lambda^{\pi}(s,a)$$=$$\sum_{k=0}^{\infty}(\mathcal{P}^{\pi})_{s_0,s}^k\pi_s(a)$$=$$\xi^{\pi}_s\pi_s(a)$ which satisfies $\lambda^{\pi}(s,a)$$\geq$$0$ and $\xi^{\pi}_s$$=$$\sum_{a\in\mathcal{A}(s)}\lambda^{\pi}(s,a)$.
Note that $\xi^{\pi}_s$$<$$\infty$ since the state $s$$\in$$S\backslash C$ is transient in $\mathcal{M}^{\pi}$. We also have
$\pi_s(a)$$=$$\frac{\lambda^{\pi}(s,a)}{\sum_{a\in\mathcal{A}(s)}\lambda^{\pi}(s,a)}$
since $\xi^{\pi}_s$$=$$\sum_{a\in\mathcal{A}(s)}\lambda^{\pi}(s,a)$$>$$0$ for reachable states. Plugging $\xi^{\pi}_s$ and $\pi_s(a)$ into \eqref{objective_MDP}, we obtain
\begin{align}\label{continuous_objective}
 &H(\mathcal{M})=\sup_{\substack{\sum_{a\in \mathcal{A}(s)}\lambda^{\pi}(s,a)> 0\\ \lambda^{\pi}(s,a)\geq 0 }}-\sum_{s\in S\backslash C}\sum_{t\in S}\nonumber\\
 &\Big[\sum_{a\in\mathcal{A}(s)}\lambda^{\pi}(s,a)\mathbb{P}_{s,a,t}\Big] \log\frac{\Big[\sum_{a'\in\mathcal{A}(s)}\lambda^{\pi}(s,a')\mathbb{P}_{s,a',t}\Big]}{\Big[\sum_{a^{''}\in\mathcal{A}(s)}\lambda^{\pi}(s,a^{''})\Big]}.
\end{align}

 Let $M$$:=$$\sup_{\pi\in\Pi(\mathcal{M})}\sum_{a\in\mathcal{A}(s)}\lambda^{\pi}(s,a)$$<$$\infty$. Then, the function $H(\mathcal{M})$ is continuous in $\lambda^{\pi}(s,a)$ and bounded over the region $\mathcal{R}$$:=$$\{ (\lambda^{\pi}(s,a))_{a\in \mathcal{A}(s)} | \lambda^{\pi}(s,a)$$\geq$$0, \sum_{a\in\mathcal{A}(s)}\lambda^{\pi}(s,a)\leq M \}$ where, if $\sum_{a\in \mathcal{A}(s)}\lambda^{\pi}(s,a)$$=$$0$, we use the convention $0\log \frac{0}{0}$$=$$0$ which preserves continuity. Note that the set $\mathcal{R}$ is closed. It is also compact since we have $\max_{\pi\in\Pi(\mathcal{M}) }\xi^{\pi}_s$$=$$\sup_{\pi\in\Pi(\mathcal{M}) }\xi^{\pi}_s$ for all $s$$\in$$S\backslash C$, which can be shown by formulating a reward maximization problem and noting that the maximum expected reward is attainable by deterministic stationary policies. We omit the explicit construction of the reward maximization problem here for brevity and refer the reader to Chapter 2 in \cite{Bertsekas} for details. Finally, since we have a continuous function maximized over a compact set in the right hand side of \eqref{continuous_objective}, the supremum is achievable. $\Box$

\noindent\textbf{Proof of Proposition \ref{modified_MDP_prop}:} Since $\mathcal{M}$ has a finite maximum entropy, all states $s$$\in$$C$ have a single successor state, i.e., $|Succ(s,a)|$$=$$1$, due to Theorem \ref{Theorem1}. Additionally, all MECs are BSC due to Theorem \ref{Theorem1}. Hence, all states $s$$\in$$C$ are either unreachable or recurrent, and have zero local entropy $L^{\pi}(s)$$=$$0$ in any MC $\mathcal{M}^{\pi}$ induced by a policy $\pi$$\in$$\Pi^S(\mathcal{M})$. Recall that for MCs with finite total entropy, the sum in \eqref{MC_entropy_biondi} is taken only over the transient states. Therefore, changing the successors of the states $s$$\in$$C$ does not affect the maximum entropy of $\mathcal{M}$ as long as the conditions $|Succ(s,a)|$$=$$1$ and $Succ(s,a)$$\subseteq$$C$ are not violated. Making states in MECs absorbing does not violate these conditions, and hence, the result follows.$\Box$\\
\noindent\textbf{Proof of Proposition \ref{prop_convex}:}  All constraints are affine in variables $\lambda(s,a)$ and $\lambda(s)$. Hence, we need only to show that the objective function is concave over the domain $\lambda(s,a)$$\geq$$0$. For a given state $s$$\in$$S$$\backslash$$C$, define the vectors $\boldsymbol{\Gamma}_s$$=$$(\eta(s,t))_{t\in S}$ and $\boldsymbol{N}_s$$=$$\nu(s)\textbf{1}$ where $\textbf{1}$$\in$$\mathbb{R}^n$. Recalling that $\eta(s,t)$ and $\nu(s)$ are functions of $\lambda(s,a)$, we define the function $f(\boldsymbol{\Gamma}_s,\boldsymbol{N}_s)$$:=$$\sum_{t\in S}\eta(s,t)\log\Big(\frac{\eta(s,t)}{\nu(s)}\Big)$
over the domain $\nu(s)$$=$$\sum_{a\in\mathcal{A}(s)}\lambda(s,a)$$\geq$$0$ and use the convention (based on continuity arguments) that $f(\boldsymbol{\Gamma}_s,\boldsymbol{0})$$=$$0$.

The function $f(\boldsymbol{\Gamma}_s,\boldsymbol{N}_s)$ is the relative entropy between the vectors $\boldsymbol{\Gamma}_s$ and $\boldsymbol{N}_s$, and thus, it is convex over the domain $\sum_{a\in\mathcal{A}(s)}\lambda(s,a)$$>$$0$ \cite{Boyd}. Moreover, since $\nu(s)$$\geq$$\eta(s,t)$$\geq$$0$ for all $s$$\in$$S\backslash C$ and $t$$\in$$S$, $f(\boldsymbol{\Gamma}_s,\boldsymbol{N}_s)$$\leq$$0$ for all $s$$\in$$S\backslash C$. Therefore, for states $s$$\in$$S\backslash C$, we can include the point $\sum_{a\in\mathcal{A}(s)}\lambda(s,a)$$=$$0$ to the domain over which the function $f$ is convex.
Now, note that the objective function in \eqref{non_reach_objective} is equal to 
$\sum_{s\in S\backslash C}-f(\boldsymbol{\Gamma}_s,\boldsymbol{N}_s)$.
Since the sum of convex functions is convex and the negation of a convex function is concave \cite{Boyd}, the objective function \eqref{non_reach_objective} is concave over the domain $\lambda(s,a)$$\geq$$0$.$\Box$\\
\noindent \textbf{Proof of Theorem \ref{Main_theorem_2}:} Assuming $H(\mathcal{M})$$<$$\infty$, we have $H(\mathcal{M}')$=$H(\mathcal{M})$ due to Proposition \ref{modified_MDP_prop}, and hence, an optimal policy for $\mathcal{M'}$ is also optimal for $\mathcal{M}$. We first prove that for a given modified MDP $\mathcal{M'}$, the objective function \eqref{non_reach_objective} of the convex program in \eqref{non_reach_objective}-\eqref{non_reach_cons6} is the maximum entropy $H(\mathcal{M'})$ of $\mathcal{M'}$. Then, we construct an optimal policy for $\mathcal{M}$ using the optimal variables $\lambda^{\star}(s,a)$ that solve the program in \eqref{non_reach_objective}-\eqref{non_reach_cons6} for ($\mathcal{M'}$, $C$), where $C$ is the set of all states in MECs in $\mathcal{M'}$. 
\par We utilize the results of \cite{Marta} to relate the variables $\lambda(s,a)$ with the expected residence time in states. 
In \cite{Marta}, it is shown that variables 
$\lambda(s,a)$$=$$\sum_{k=0}^{\infty}(\mathcal{P}^{\pi})_{s_0,s}^k\pi_s(a)$$=$$\xi^{\pi}_s\pi_s(a)$
satisfy the constraint in \eqref{non_reach_cons1} and corresponds to the expected residence time in a state-action pair $(s,a)$ in an induced MC $\mathcal{M}^{\pi}$. Additionally, $\lambda(s)$ corresponds to the reachability probability of states $s$$\in$$C$.
Then, it is clear that for states $s$$\in$$S\backslash C$, 
{\setlength{\mathindent}{-2pt}
\begin{align}
\label{state_residence}
\xi^{\pi}_s=\sum_{a\in\mathcal{A}(s)}\lambda(s,a).
\end{align}}\noindent
Additionally, if $\sum_{a\in\mathcal{A}(s)}\lambda(s,a)$$>$$0$, we have
{\setlength{\mathindent}{-2pt}
\begin{align}\label{policy_construction}
\pi_s(a)=\frac{\lambda(s,a)}{\sum_{a\in\mathcal{A}(s)}\lambda(s,a)}.
\end{align}}
\par Recall that for all $s$$\in$$C$ and $\pi$$\in$$\Pi^S(\mathcal{M}')$, we have $L^{\pi}(s)$$=$$0$ since $H(\mathcal{M}')$$<$$\infty$. Therefore, 
\begin{align}
H(\mathcal{M}')&=\sup_{\pi\in{\Pi^S(\mathcal{M}')}}\Big[\sum_{s\in S}\xi^{\pi}_sL^{\pi}(s)\Big]\\
\label{equivalence_in_objective}
&=\max_{\pi\in{\Pi^S(\mathcal{M}')}}\Big[\sum_{s\in S\backslash C}\xi^{\pi}_sL^{\pi}(s)\Big],
\end{align}
\noindent due to Proposition \ref{Biondi_theorem} and Theorem \ref{Theorem1}. Our aim is to show that the expression in $\eqref{equivalence_in_objective}$ is equal to the objective in \eqref{non_reach_objective}.
\par For an arbitrary $\pi$$\in$$\Pi^S(\mathcal{M}')$, define the set $G^{\pi}$$:=$$\{s$$\in$$ S\backslash C | \xi^{\pi}_s$$=$$0 \}$. Note that if $G^{\pi}$$\neq$$\emptyset$ for some $\pi$$\in$$\Pi^S(\mathcal{M}')$, the states $s$$\in$$G^{\pi}$ do not affect the value of \eqref{equivalence_in_objective} by the definition of $G^{\pi}$. We consider two cases: (1) $G^{\pi}$$\neq$$\emptyset$ and (2) $G^{\pi}$$=$$\emptyset$. For case 1, we will show that states $s$$\in$$G^{\pi}$ do not affect the value of \eqref{non_reach_objective}, and for case 2, we will show that the expression in \eqref{equivalence_in_objective} is equal to the objective in \eqref{non_reach_objective}.

\textbf{Case 1}: We assume that $G^{\pi}$$\neq$$\emptyset$ and show that for any $s$$\in$$G^{\pi}$,
\begin{align}
\label{limit_first_step}
 \sum_{t\in S}\eta(s,t)\log\Big(\frac{\eta(s,t)}{\nu(s)}\Big)=0.
\end{align}
\par Considering \eqref{non_reach_cons3}-\eqref{non_reach_cons5}, and noting that $0$$\leq$$\mathbb{P}_{s,a,t}$$\leq$$1$ for all $t$$\in$$S$, we have $\nu(s)$$\geq$$\eta(s,t)$$\geq$$0$ for all $t$$\in$$S$. Hence, for any $s$$\in$$G^{\pi}$, we have $\nu(s)$$=$$\eta(s,t)$$=$$0$ for all $t$$\in$$S$ due to the definition of the set $G^{\pi}$ and \eqref{state_residence}. 
We use the convention $0\log\frac{0}{0}$$=$$0$ which is based on continuity arguments and the fact that whenever $\nu(s)$$=$$0$, we have $\eta(s,t)$$=$$0$ for all $t$$\in$$S$. As a result, we conclude that the states $s$$\in$$G^{\pi}$ do not affect the value of the objective in \eqref{non_reach_objective}.

\textbf{Case 2}: We assume that $G^{\pi}$$=$$\emptyset$. In this case, for any $\pi$$\in$$\Pi^S(\mathcal{M}')$, we have $\xi^{\pi}_s$$=$$\sum_{a\in\mathcal{A}(s)}\lambda(s,a)$$>$$0$ and, \eqref{policy_construction} holds for all $s$$\in$$S\backslash$$C$ and $a$$\in$$\mathcal{A}(s)$. Plugging \eqref{state_residence} and \eqref{policy_construction} into \eqref{equivalence_in_objective}, we obtain the objective function in \eqref{non_reach_objective}. (Note that $\eta(s,t)$ and $\nu(s)$ variables can be written in terms of $\lambda(s,a)$ using \eqref{non_reach_cons3}-\eqref{non_reach_cons4}.) We conclude that the problem in \eqref{non_reach_objective}-\eqref{non_reach_cons6} computes the maximum entropy of $\mathcal{M}'$.
\par Now, we construct an optimal policy for $\mathcal{M}$. We show in \eqref{equivalence_in_objective} that states $s$$\in$$C$ does not affect the value of $H(\mathcal{M'})$. Therefore, an arbitrary assignment of actions in states $s$$\in$$C$ does not affect the optimality of a policy. Similarly, for a given optimal policy $\pi^{\star}$, an arbitrary assignment of actions in states $s$$\in$$G^{\pi^{\star}}$ does not affect the optimality due to \eqref{limit_first_step}. Additionally, by the construction given in \eqref{policy_construction}, an optimal policy for states $s$$\in$$S\backslash(C$$\cup$$ G^{\pi^{\star}})$ satisfies 
$\pi^{\star}_s(a)$$=$$\frac{\lambda^{\star}(s,a)}{\sum_{a\in\mathcal{A}(s)}\lambda^{\star}(s,a)}$,
where $\lambda^{\star}(s,a)$ are optimal variables for the problem in \eqref{non_reach_objective}-\eqref{non_reach_cons6}. Since an optimal policy for $\mathcal{M}'$ is also optimal for $\mathcal{M}$ due to Proposition \ref{modified_MDP_prop}, we conclude that Algorithm \ref{Algo_2} returns an optimal policy for $\mathcal{M}$. $\Box$\\
\noindent\textbf{Proof of Lemma \ref{paths_lemma}:}  
For an MC $\mathcal{M}^{\pi}$ induced by a policy $\pi$$\in$$\Pi^S(\mathcal{M})$, let $S_B$ and $S_{B_0}$ be the union of its BSCCs and the set of its transient states, respectively. Moreover, let $T$$:=$$Paths^{\pi}_{fin}(\mathcal{M})$$\cap$$ (S\backslash S_B)^{\star}S_B$. For states $s$$\in$$S_{B_0}$ and $t$$\in$$S$, define sets
\vspace{-0.4cm}
\begin{align*}
&A_{s,k}:=\{s_0\ldots s_n \in T :\ n\in\mathbb{N},\  \sum_{i=0}^n\mathbbm{1}_{\{s_i=s\}}=k \},\\
&B_{(s,t),k}:=\{s_0\ldots s_n \in T : \ n\in\mathbb{N}, \sum_{i=1}^n\mathbbm{1}_{\{(s_{i-1},s_{i})=(s,t)\}}=k \}.
\end{align*}
\vspace{-0.4cm}

Note that $A_{s,k}$ is the collection of all paths along which the state $s$ is $k$ times visited and a BSCC in $\mathcal{M}^{\pi}$ is reached. Similarly, the set $B_{(s,t),k}$ is the collection of all paths along which the edge between state $s$ and state $t$ is $k$ times traversed, and a BSCC in $\mathcal{M}^{\pi}$ is reached. 

It is known that any finite MC almost surely reaches a BSCC \cite{Model_checking}. Thus, for any $s$$\in$$S_{B_0}$ and $t$$\in$$S$, we have
\vspace{-0.2cm}
\begin{align}\label{summation_one}
\sum_{k=0}^{\infty}\text{Pr}^{\pi}_{\mathcal{M}}(A_{s,k})=\sum_{k=0}^{\infty}\text{Pr}^{\pi}_{\mathcal{M}}(B_{(s,t),k})=1.
\end{align}

One can show using \eqref{summation_one} that
\vspace{-0.2cm}
\begin{align}
\label{state_res}
\xi^{\pi}_s=\sum_{k=0}^{\infty} k\ \text{Pr}^{\pi}_{\mathcal{M}}(A_{s,k})
\end{align}
\vspace{-0.2cm}

for transient states $s$$\in$$S_{B_0}$ in $\mathcal{M}^{\pi}$. (We omit the derivation here. The result can be obtained by using the countability of $A_{s,k}$ and performing a series of algebraic manipulations to obtain \eqref{residence}. A similar derivation can also be found in \cite{Saerens}.)

Let $\xi^{\pi}_{s,t}$ denote the expected number of transitions from a state $s$$\in$$S_{B_0}$ to state $t$$\in$$S$. Then, we have

\vspace{-0.4cm}
\begin{align}
\label{transition_res}
\xi^{\pi}_{s,t}=\sum_{k=0}^{\infty} k\ \text{Pr}^{\pi}_{\mathcal{M}}(B_{(s,t),k}),
\end{align} 
analogously to \eqref{state_res}. Additionally, the relation between \eqref{state_res} and \eqref{transition_res} is given by
$\mathcal{P}^{\pi}_{s,t}\xi^{\pi}_{s}=\xi^{\pi}_{s,t}$, which can be obtained by using a method similar to the one used in \cite{Serfozo} to derive \eqref{residence}. 

Let $N^{s_0\ldots s_n}_{s,t}$ be the number of transitions made from state $s$$\in$$S_{B_0}$ to state $t$$\in$$S$ along a finite path fragment $s_0\ldots s_n$$\in$$T$. Then, we have
\begin{align}\label{algebra_property}
\mathcal{P}^{\pi}_{s,t}\xi^{\pi}_{s}=\xi^{\pi}_{s,t}&=\sum_{k=0}^{\infty}\sum_{s_0\ldots s_n\in B_{(s,t),k}} k\  \mathcal{P}^{\pi}(s_0\ldots s_n)\\
\label{residence_equivalence}
&=\sum_{s_0\ldots s_n\in T} N^{s_0\ldots s_n}_{s,t}\  \mathcal{P}^{\pi}(s_0\ldots s_n),
\end{align} 
where the equality in \eqref{algebra_property} follows from the fact that set $B_{(s,t),k}$ is countable and each element $s_0\ldots s_n$$\in$$B_{(s,t),k}$ is measurable. The equality in \eqref{residence_equivalence} is due to the fact that any finite path fragment $s_0\ldots s_n\in T$ is an element of one and only one set $B_{(s,t),k}$, and that for a given path fragment $s_0\ldots s_n$, we have $k$$=$$N^{s_0\ldots s_n}_{s,t}$ by definition. 

We next express the probability of a finite path fragment in terms of the number of transition appearances.  Then, we have
\begin{align}
\label{measure_2}
\mathcal{P}^{\pi}(s_0\ldots s_m)=\prod_{(s,t)\in S_{B_0}\times S : N^{s_0\ldots s_m}_{s,t}>0}(\mathcal{P}^{\pi}_{s,t})^{N^{s_0\ldots s_m}_{s,t}}.
\end{align}

By assumption, we have $H(\mathcal{M},{\pi})$$<$$\infty$. If $s_0$$\in$$S_B$, both the entropy and the entropy of paths for $\mathcal{M}^{\pi}$ are zero; hence, we only analyze the case $s_0$$\in$$S_{B_0}$. In this case, the summation in \eqref{MC_entropy_biondi} is taken over transient states $s$$\in$$S_{B_0}$ since $H(\mathcal{M},{\pi})$$<$$\infty$. As a result, 
\begin{align}\label{entropy_def_derive}
H(\mathcal{M},{\pi})&=-\sum_{s\in S_{B_0}}\xi^{\pi}_s\sum_{t\in S}\mathcal{P}^{\pi}_{s,t}\log \mathcal{P}^{\pi}_{s,t}\\ 
\label{remove_non_trans}
&=-\sum_{(s,t)\in S_{B_0}\times S : \mathcal{P}^{\pi}_{s,t}>0}\xi^{\pi}_s\mathcal{P}^{\pi}_{s,t}\log \mathcal{P}^{\pi}_{s,t}\\
\label{trans_res_derive}
&=-\sum_{(s,t)\in S_{B_0}\times S : \mathcal{P}^{\pi}_{s,t}>0}\nonumber \\
&\qquad \sum_{s_0\ldots s_n\in T}N^{s_0\ldots s_n}_{s,t}\ \mathcal{P}^{\pi}(s_0\ldots s_n)\log \mathcal{P}^{\pi}_{s,t}\\
\label{remove_zeros}
&= -\sum_{(s,t)\in S_{B_0}\times S : N^{s_0\ldots s_n}_{s,t}>0}\nonumber\\
&\qquad \sum_{s_0\ldots s_n\in T}N^{s_0\ldots s_n}_{s,t}\ \mathcal{P}^{\pi}(s_0\ldots s_n)\log \mathcal{P}^{\pi}_{s,t},
\end{align}
where \eqref{remove_non_trans} follows by removing transitions $\mathcal{P}^{\pi}_{s,t}$$=$$0$ and using the convention $ 0\log 0$$=$$0$, \eqref{trans_res_derive} follows from \eqref{residence_equivalence}, and \eqref{remove_zeros} is obtained by removing state pairs $(s,t)$$\in$$S_{B_0}$$\times$$S$ for which $N^{s_0\ldots s_n}_{s,t}$$=$$0$. 

Now, we analyze the entropy of paths. The entropy of paths $H(Paths^{\pi}(\mathcal{M}))$ for the induced MC $\mathcal{M}^{\pi}$ can be written as

\begin{align}
\label{first_use_measure}
&H(Paths^{\pi}(\mathcal{M}))=\nonumber\\
&-\sum_{s_0\ldots s_n\in T}\mathcal{P}^{\pi}(s_0\ldots s_n)\log \mathcal{P}^{\pi}(s_0\ldots s_n),\\
\label{use_measure}
&=-\sum_{s_0\ldots s_m\in T}\nonumber \\
&\sum_{(s,t)\in S_{B_0}\times S : N^{s_0\ldots s_n}_{s,t}>0}N^{s_0\ldots s_n}_{s,t}\mathcal{P}^{\pi}(s_0\ldots s_n)\log\mathcal{P}^{\pi}_{s,t},
\end{align}
where \eqref{use_measure} is obtained by plugging \eqref{measure_2} into \eqref{first_use_measure}. \\
Since \eqref{remove_zeros} and \eqref{use_measure} are equal, we conclude that $H(\mathcal{M},{\pi})$$=$$H(Paths^{\pi}(\mathcal{M}))$ for any $\pi$$\in$$\Pi^S(\mathcal{M})$, under the assumption that $H(\mathcal{M},{\pi})$$<$$\infty$ for all $\pi$$\in$$\Pi^S(\mathcal{M})$. $\Box$


\section{}\label{infinite_case_appendix}

In this appendix, we provide procedures to solve entropy maximization and constrained entropy maximization problems for MDPs with infinite maximum entropy. 
\subsubsection{Entropy maximization}\label{entropy max sec}
In this case, for a given MDP $\mathcal{M}$ with the union $(C, D)$ of its MECs, there exists at least one state $s^{\star}$$\in$$C$ such that $\lvert \cup_{a\in D(s^{\star})} Succ(s^{\star},a)\rvert$$>$$1$ due to Theorem \ref{Theorem1}. We aim to synthesize a policy that induces an MC where the state $s^{\star}$ is both stochastic and recurrent. For simplicity, we assume that  there exists only one state $s^{\star}$ such that  $\lvert \cup_{a\in D(s^{\star})} Succ(s^{\star},a)\rvert$$>$$1$. We form the modified MDP $\mathcal{M}'$ by replacing each BSC MEC in $\mathcal{M}$ with an absorbing state. Let $S_B$ and $S_{NB}$ be the set of all states in BSC MECs and non-BSC MECs in $\mathcal{M}$, respectively. We consider two cases, namely $s^{\star}$$\in$$S_B$ and $s^{\star}$$\in$$S_{NB}$. If $s^{\star}$$\in$$S_B$, let $C'$ be the union of all absorbing states in $\mathcal{M}'$ that are replaced with BSC MECs in $\mathcal{M}$, and $C_{\star}$ be the absorbing state that is replaced with the MEC that $s^{\star}$ is contained in. We solve the problem in \eqref{non_reach_objective}-\eqref{non_reach_cons6} for ($\mathcal{M}'$, $C'$) together with the constraint $\lambda(C_{\star})$$>$$0$. (Note that if there is a non-BSC MEC in $\mathcal{M}'$, the constraint \eqref{residence_bound} should also be included to this optimization problem.) We then use step 3 of Algorithm \ref{Algo_2} to obtain a policy for states $s$$\not\in$$C'$, and choose actions in state $s^{\star}$ such that $\lvert Succ(s^{\star})\rvert$$>$$1$ in the induced MC. By construction, the state $s^{\star}$ is both stochastic and recurrent in the induced MC, and due to Proposition \ref{Biondi_theorem}, the entropy of the induced MC is infinite. If $s^{\star}$$\in$$S_{NB}$, we replace the MEC that state $s^{\star}$ is contained in with an absorbing state  and follow steps similar to the ones in the case $s^{\star}$$\in$$S_{B}$ to obtain an optimal policy.

\subsubsection{Constrained entropy maximization}
We suppose that the feasible policy space for the problem in \eqref{max_ent_product_objective}-\eqref{LTL_constraint_product} is not empty. The product MDP $\mathcal{M}_p$ contains a MEC $(C,D)$ such that  $|\cup_{a\in D(s^{\star})}Succ(s^{\star},a)|$$>$$1$ for some $s^{\star}$$\in$$C$ due to Theorem \ref{Theorem1}. We assume that there exists at least one non-BSC MEC in $\mathcal{M}_p$ and there is only one state $s^{\star}$ in $\mathcal{M}_p$ such that $|\cup_{a\in D(s^{\star})}Succ(s^{\star},a)|$$>$$1$. These assumptions are introduced just to simplify the case analysis. We first partition the states into three disjoint sets $B,S_0$ and $S_r$ as explained in Section \ref{product_policy_section}. Then, we form the modified MDP $\mathcal{M}_p'$ by replacing each BSC MEC in $\mathcal{M}_p$ with an absorbing state. Let $S_B$ and $S_{NB}$ be the set of all states in BSC MECs and non-BSC MECs in $\mathcal{M}_p$, respectively. We consider two cases: 1) $s^{\star}$$\in$$S_B$ and 2) $s^{\star}$$\in$$S_{NB}$.

\textit{Case 1:}  
 If $s^{\star}$$\in$$S_B$, let $C'$ be the union of all absorbing states in $\mathcal{M}_p'$ that are replaced with BSC MECs in $\mathcal{M}_p$, and $C_{\star}$ be the absorbing state that is replaced with the MEC that $s^{\star}$ is contained in. We obtain a policy for states $s$$\in$$S_r$ by solving the problem in \eqref{non_reach_objective}-\eqref{non_reach_cons6} for ($\mathcal{M}_p'$,$C'$,$\beta$,$\Gamma$) together with the constraints \eqref{residence_bound}, \eqref{reach_cons} and $\lambda(C_{\star})$$>$0. Then, we choose actions in state $s^{\star}$ such that $\rvert Succ(s^{\star})\lvert$$>$$1$ in the induced MC. Note that if this problem is infeasible, then there exists no policy that induces an MC with infinite entropy whose paths satisfies the LTL specification with probability $\beta$. In this case, the maximum constrained entropy is unbounded, and we follow the steps that are explained in Section \ref{unbounded_constrained_case} to synthesize a policy that induces an MC with arbitrarily large entropy. 

\textit{Case 2:} If $s^{\star}$$\in$$S_{NB}$, we consider two cases, namely $s^{\star}$$\in$$B$$\cup$$S_0$ and $s^{\star}$$\in$$S_r$. If $s^{\star}$$\in$$B\cup S_0$, we replace the MEC that $s^{\star}$ belongs to in $\mathcal{M}_p'$ with an absorbing state. Then, we synthesize a policy that induces an MC with infinite entropy whose paths satisfy the LTL specification with probability $\beta$ as explained in Case 1. Additionally, to ensure that the state $s^{\star}$ is recurrent in the induced MC, we choose actions in states that belong to the same MEC with $s^{\star}$ such that the MEC forms a BSCC in the induced MC. If $s^{\star}$$\in$$S_r$, the maximum constrained entropy is not infinite because no state $s$$\in$$S_r$ can be recurrent in an induced MC that satisfies the LTL specification with probability $\beta$$>$$0$. In this case, the maximum constrained entropy is unbounded, and we use the procedure explained in Section \ref{unbounded_constrained_case} to synthesize a policy that induces an MC with arbitrarily large entropy.

\vspace{-1cm}
\begin{IEEEbiography}[{\includegraphics[width=1in,height=1.25in,clip,keepaspectratio]{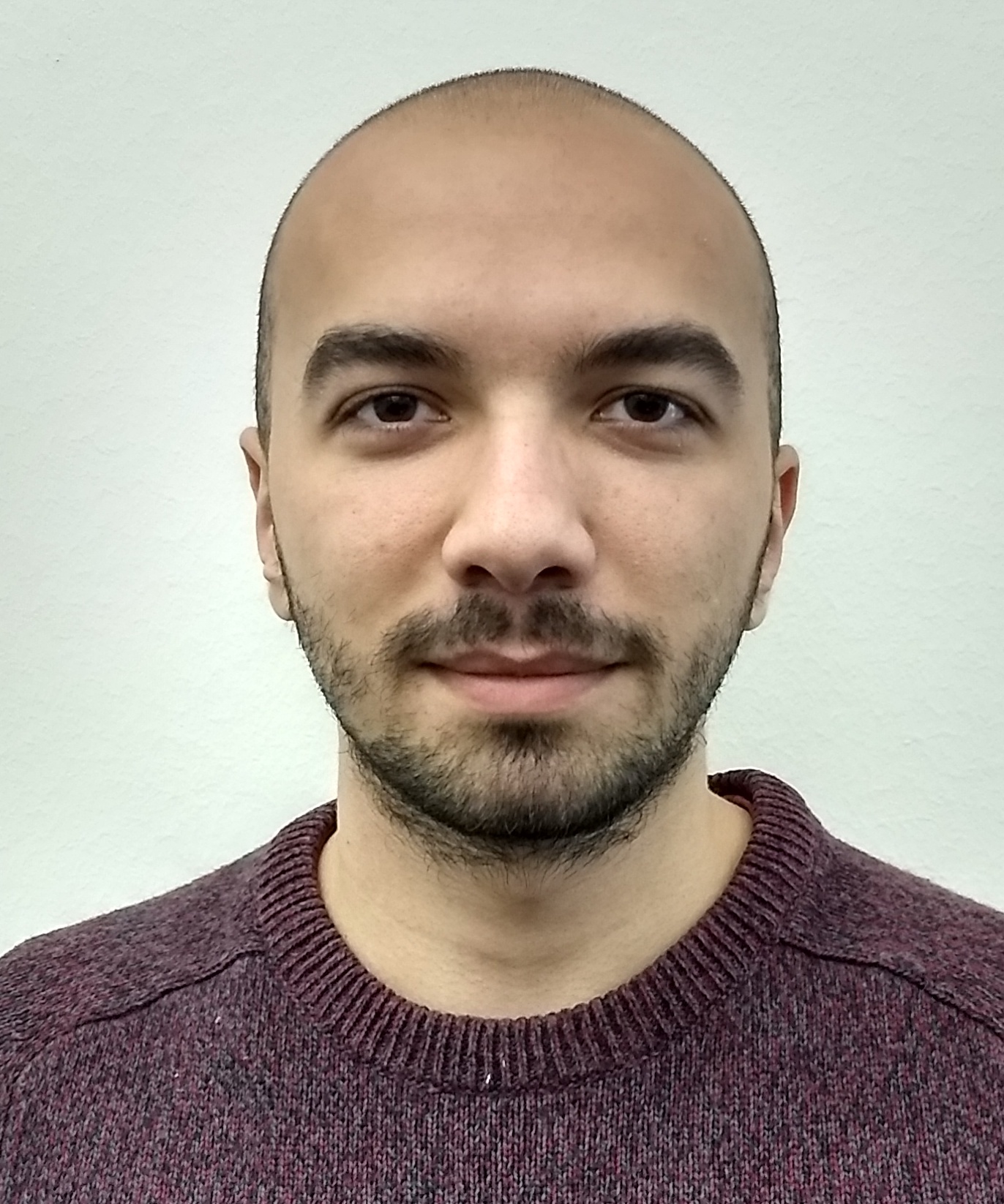}}]{Yagiz Savas} joined the Department of Aerospace Engineering at the University of Texas at Austin as a Ph.D. student in Fall 2017. He received his B.S. degree in Mechanical Engineering from Bogazici University in 2017. His research focuses on developing theory and algorithms that guarantee desirable behavior of autonomous systems operating in adversarial environments.
\end{IEEEbiography}
\vspace{-1.2cm}
\begin{IEEEbiography}[{\includegraphics[width=1in,height=1.25in,clip,keepaspectratio]{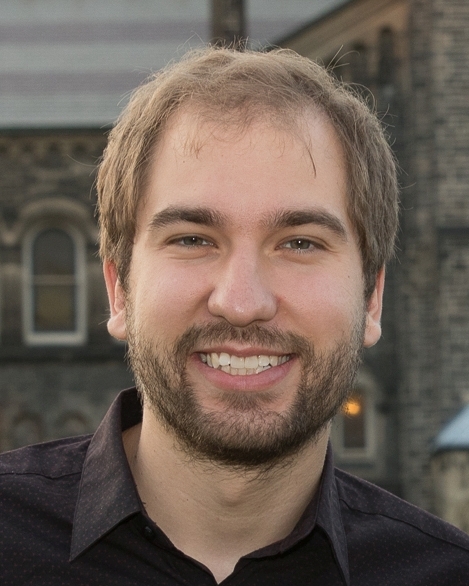}}] {Melkior Ornik} is an assistant professor in the Department of Aerospace Engineering and the Coordinated Science Laboratory at the University of Illinois at Urbana-Champaign. He received his Ph.D. degree from the University of Toronto in 2017. His research focuses on developing theory and algorithms for learning and planning of autonomous systems operating in uncertain, complex and changing environments, as well as in scenarios where only limited knowledge of the system is available.
\end{IEEEbiography}
\vspace{-1.4cm}
\begin{IEEEbiography}[{\includegraphics[width=1in,height=1.25in,clip,keepaspectratio]{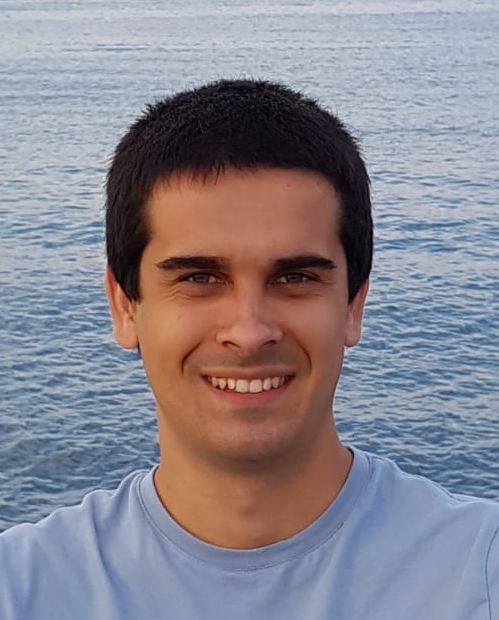}}]{Murat Cubuktepe} joined the Department of Aerospace Engineering at the University of Texas at Austin as a Ph.D. student in Fall 2015. He received his B.S degree in Mechanical Engineering from Bogazici University in 2015. His main current research interests are verification and synthesis of uncertain, parametric and partially observable probabilistic systems. He also focuses on applications of convex optimization in formal methods and controls.
\end{IEEEbiography}
\vspace{-1.2cm}
\begin{IEEEbiography}[{\includegraphics[width=1in,height=1.25in,clip,keepaspectratio]{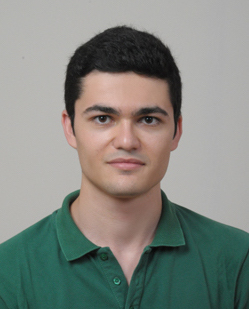}}]{Mustafa O. Karabag} joined the Department of Electrical and Computer Engineering at the University of Texas at Austin as a Ph.D. student in Fall 2017. He received his B.S. degree in Electrical and Electronics Engineering from Bogazici University in 2017. His research focuses on developing theory and algorithms for non-inferable planning in adversarial environments.
\end{IEEEbiography}
\vspace{-1.2cm}
\begin{IEEEbiography}[{\includegraphics[width=1in,height=1.25in,clip,keepaspectratio]{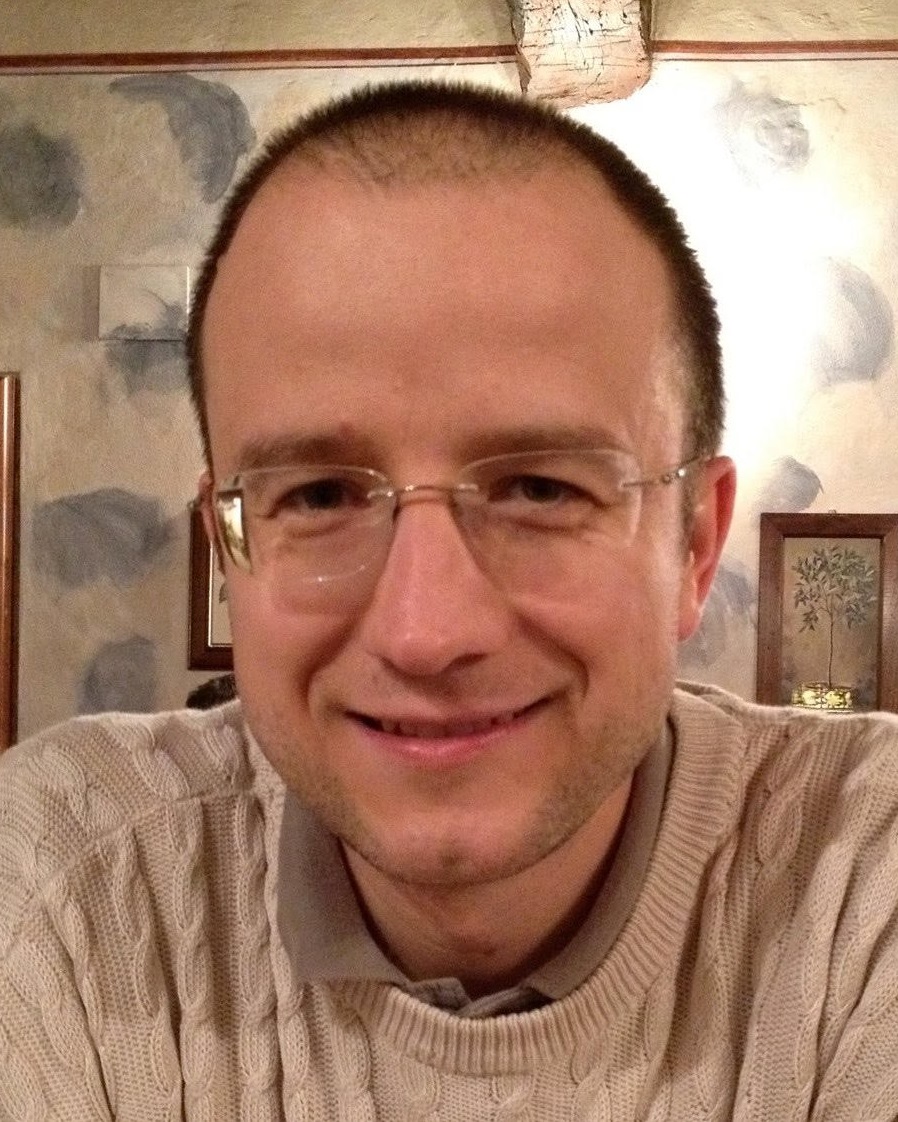}}]{Ufuk Topcu} joined the Department of Aerospace Engineering at the University of Texas at Austin as an assistant professor in Fall 2015. He received his Ph.D. degree from the University of California at Berkeley in 2008. He held research positions at the University of Pennsylvania and California Institute of Technology. His research focuses on the theoretical, algorithmic and computational aspects of design and verification of autonomous systems through novel connections between formal methods, learning theory and controls.
\end{IEEEbiography}

 \end{document}